\documentclass[11pt]{refart}
   \settextfraction{0.75}
   \newenvironment{SCfigure*}{\begin{figure*}}{\end{figure*}}

\let\LTXsect\section
\renewcommand{\section}{\clearpage\LTXsect}

\title{Multiscale modelling of microscale heterogeneous systems: 
analysis supports systematic and efficient macroscale
modelling and simulation}
\author{A.J. Roberts
\thanks{Supported by Australian Research Council grants DP180100050, DP150102385, et al.}\\School of Mathematical Sciences, University of Adelaide
\\\protect\url{mailto:anthony.roberts@adelaide.edu.au}
}
\date{\today}

\usepackage[nohyperref]{ajr}
\usepackage{microtype,amsmath,amsthm,defns,MnSymbol,enumitem,multicol,needspace}
\setlist{parsep=0.5ex}
\usepackage{pgfplots}\pgfplotsset{compat=newest}
\usepackage{mybiblatex,mycleveref}
\usepackage[twodeep]{etocx}

\let\matlab\relax
\usepackage{fancyvrb}
\newenvironment{matlab}%
    {\Verbatim[numbers=left,xleftmargin=1em]}%
    {\endVerbatim}

\Vec f \Vec u \Vec v \Vec x \Vec y
\Vec U
\def\ad{\hat a}
\def\am{\bar a}
\newcounter{i}
\def\fx{x}
\def\fu{\mathfrak u}
\Bb X
\newcommand{\E}[2]{$#1\cdot10^{#2}$}
\def\script{\Matlab\slash Octave}

\hypersetup{colorlinks
,linkcolor=RoyalBlue,citecolor=RoyalBlue,pagecolor=RoyalBlue,%
            urlcolor=magenta,filecolor=magenta}%,breaklinks,%
\begin{document}

\maketitle

\begin{abstract}
These are lecture notes for five sessions in the AMSI Winter School on \emph{Computational Modelling of Heterogeneous Media} held at QUT in July 2019 [\url{https://ws.amsi.org.au/}].

Aim: Discuss a mix of new mathematical approaches for multiscale modelling, heterogeneous material in particular, along with corresponding novel computational techniques and issues. I include discussion of a developing toolbox that empowers you to implement effective multiscale `equation-free' computation 
\\{[\url{https://github.com/uoa1184615/EquationFreeGit.git}].}
\end{abstract}

\paragraph{Pre-requisites} undergraduate ordinary and partial differential equations, state space, trajectories, stability, bifurcations, power series solutions, separation of variables, eigen-problems and Sturm--Liouville theory; basic numerical methods for time integration of ODEs and for spatial discretisation PDEs, some perturbation methodology. 

\paragraph{Pre-reading} recommend the review article by \cite{Kevrekidis09a}.  Perhaps also get an idea of centre manifolds 
\\{[\url{https://en.wikipedia.org/wiki/Center_manifold}].}

\paragraph{Activities} sprinkled are some small activities for you. 

\tableofcontents

\section{Introducing a powerful approach to macroscale modelling} 
\label{sammt}
\localtableofcontents

Let's revisit homogenisation in some simple problems in order to clarify and resolve macroscale modelling issues.
Here we primarily discuss diffusion across an inhomogeneous lattice, and comment on cognate elastic vibrations.
Chapter~7 of my book \cite[]{Roberts2014a} discusses more details, insights and applications.

This section introduces a powerful new framework for understanding and creating multiscale computational algorithms:
\begin{itemize}
\item resolves how real~\(x\) relates to discrete lattice/cells;
\item applies to \emph{finite} microscales in a \emph{finite} macroscale;
\item finds sharp lower bounds for allowable space and time scales;
\item derives accurate boundary conditions for \pde\ models;
\item and could illuminate initial conditions, uncertainty, forcing, etc.
\end{itemize}

\subsection{Homogenise period-two diffusion on a lattice---the leading approximation}

\begin{figure*}
\caption{schematic diagram of inhomogeneous diffusion on a microscale lattice, spacing~\(d\), of material~$u_j(t)$, with two-periodic diffusivity $a$~and~$b$ between even and odd lattice points.}
\label{fig:svhld}
\centering
\setlength{\unitlength}{1ex}
\begin{picture}(75,8)
%\put(0,0){\framebox(75,8){}}
\put(0,-2){\thicklines
\put(0,5){\vector(1,0){72}$\ j$}
\multiput(3,4)(8,0){9}{\line(0,1)2}
\setcounter{i}{0}
\multiput(10,2)(8,0){7}{$\arabic{i}$\stepcounter{i}}
\multiput(2,2)(64,0){2}{$\cdots$}
\setcounter{i}{0}
\multiput(8.5,8)(8,0){7}{$u_{\arabic{i}}(t)$\stepcounter{i}}
\multiput(2,8)(64,0){2}{$\cdots$}
\multiput(6,6)(16,0)4{$a$}
\multiput(14,6)(16,0)4{$b$}
}
\end{picture}
\end{figure*}
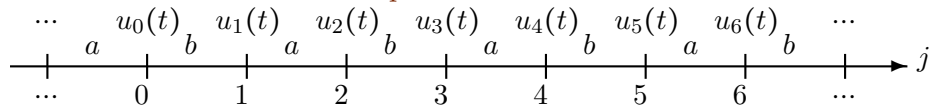%
Suppose a microscale lattice, \(x_j=jd\) as in \cref{fig:svhld}, has a property~\(u_j(t)\) that diffuses with the period-two coefficients, named~\(a,b\):
\begin{equation}
\dot u_j=\begin{cases}
a[u_{j+1}-u_j]+b[u_{j-1}-u_j]&\text{odd }j,\\
b[u_{j+1}-u_j]+a[u_{j-1}-u_j]&\text{even }j.
\end{cases}
\label{eqsvhdode}
\end{equation}

Define~\(U(x,t)\) to be an average of nearby~\(u_{\text{odd}}(t)\) and~\(u_{\text{even}}(t)\) then many arguments derive the macroscale homogenised diffusion
\begin{equation}
\D tU\approx D\DD xU\,,
\quad D=\frac{2ab}{a+b}=\frac1{(1/a+1/b)/2}.
\end{equation}
In this \pde\ model there is no difference between diffusions~\(a,b,a,\ldots\) and \(b,a,b,\ldots\) so we could imperceptibly swap---as we do in the next subsection.

The common approach of letting the microscale spacing \(d\to0\) appears a good way to derive such a \pde\ because ``\(d\to0\)'' squashes a lot of difficult issues.  
Unfortunately, these difficult issues are often of practical importance, such as determining boundary conditions for a \pde\ model (\cref{ssBCs}).

\subsection{Resolve the microscale to regularise, rather than singularise}
\label{ssrmrrs}

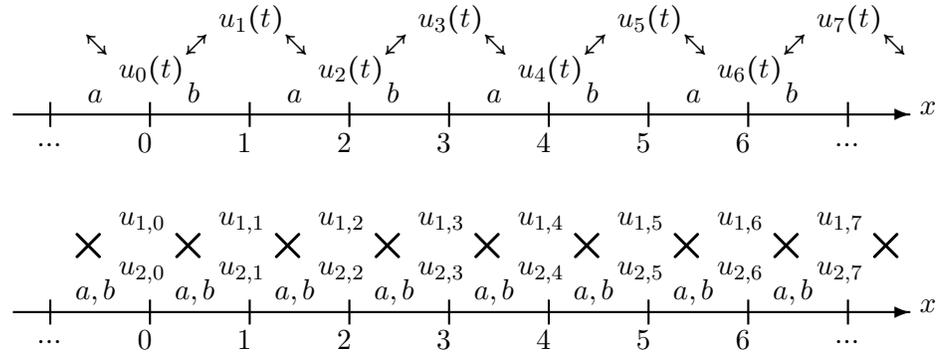
\begin{figure*}
\caption{\label{fig:dbl}First~(top) redraw \cref{fig:svhld} so that the \emph{odd} points are in the upper line and the \emph{even} points in the lower line.  Then the diffusivity~\(a\) always acts in `NW--SE' direction and the diffusivity~\(b\) always acts in `NE--SW' direction.
Second~(bottom) double the number of variables by filling in the gaps: \(u_{1,j}(t)\)~in the upper line; and \(u_{2,j}\)~in the lower line.  Here the diffusivities~\(a\) and~\(b\) always act in the same direction.  This forms a system homogeneous on the lattice!}
\centering
\setlength{\unitlength}{1ex}
\begin{tabular}{c}
\begin{picture}(75,12)
%\put(0,0){\framebox(75,12){}}
\put(0,-2){\thicklines
\put(0,5){\vector(1,0){72}$\ x$}
\multiput(3,4)(8,0){9}{\line(0,1)2}
\setcounter{i}{0}
\multiput(10,2)(8,0){7}{$\arabic{i}$\stepcounter{i}}
\multiput(2,2)(64,0){2}{$\cdots$}
\setcounter{i}{0}
\multiput(8.5,8)(16,0){4}{$u_{\arabic{i}}(t)$\stepcounter{i}\stepcounter{i}}
\setcounter{i}{1}
\multiput(16.5,12)(16,0){4}{$u_{\arabic{i}}(t)$\stepcounter{i}\stepcounter{i}}
\multiput(13.5,10)(16,0){4}{\large$\neswarrow$}
\multiput(5.5,10)(16,0){5}{\large$\nwsearrow$}
\multiput(6,6)(16,0)4{$a$}
\multiput(14,6)(16,0)4{$b$}
}
\end{picture}
\\
\begin{picture}(75,15)
%\put(0,0){\framebox(75,12){}}
\put(0,-2){\thicklines
\put(0,5){\vector(1,0){72}$\ x$}
\multiput(3,4)(8,0){9}{\line(0,1)2}
\setcounter{i}{0}
\multiput(10,2)(8,0){7}{$\arabic{i}$\stepcounter{i}}
\multiput(2,2)(64,0){2}{$\cdots$}
\setcounter{i}{0}
\multiput(8.5,8)(8,0){8}{$u_{2,\arabic{i}}$\stepcounter{i}}
\setcounter{i}{0}
\multiput(8.5,12)(8,0){8}{$u_{1,\arabic{i}}$\stepcounter{i}}
\multiput(4.5,9)(8,0){9}{\Huge$\times$}
\multiput(5,6)(8,0)8{$a,b$}
}
\end{picture}
\end{tabular}
\end{figure*}%
Let's form a system that is homogeneous on the lattice.
\cref{fig:dbl} shows how to do this by forming two decoupled systems that we analyse as one whole system.
Usefully, the combined system is homogeneous on the lattice: every lattice point has the same equations because diffusivity~\(a\) always acts in the same direction, and the diffusivity~\(b\) always acts in the other direction.
The \ode{}s are
\begin{align*}&
\dot u_{1,j}=a[u_{2,j+1}-u_{1,j}]+b[u_{2,j-1}-u_{1,j}],
\\&
\dot u_{2,j}=b[u_{1,j+1}-u_{2,j}]+a[u_{1,j-1}-u_{2,j}].
\label{eq:}
\end{align*}

Non-dimensionalise on the microscale length so that in effect \(d=1\). 
Writing the \ode{}s as in terms of the spatial location of the lattice points, \(x=j\), the differential equations for \(\uv=(u_1,u_2)\) are 
\begin{subequations}\label{eqssvhdmic}%
\begin{align}&
\dot u_1(x)=a[u_2(x+1)-u_1(x)]
\nonumber\\&\qquad{}
+b[u_2(x-1)-u_1(x)],
\\&
\dot u_2(x)=b[u_1(x+1)-u_2(x)]
\nonumber\\&\qquad{}
+a[u_1(x-1)-u_2(x)].
\end{align}
\end{subequations}
We analyse these near-rigorously, discover the homogenisation, and more.

\subsection{An ensemble of phase shifts underpins rigorous homogenisation,} 
\label{ssepsurh}

Fourier space provides a rigorous route \cite[]{Chen2014}.
However, let's employ a corresponding direct approach \cite[based upon][which also compares with Fourier transform]{Roberts2013a}.
We focus on `slowly varying' solutions by expanding in the notionally small~\(\partial_x\):  
\begin{align*}&
\dot u_1=(a+b)(u_2-u_1)+(a-b)u_{2x}
\\&\qquad{}
+\tfrac12(a+b)u_{2xx}+\cdots,
\\&
\dot u_2=(a+b)(u_1-u_2)-(a-b)u_{1x}
\\&\qquad{}
+\tfrac12(a+b)u_{1xx}+\cdots.
\end{align*}
\begin{enumerate}
\def\smi#1{\marginpar{
\begin{tikzpicture}
\begin{axis}[footnotesize,font=\footnotesize
    ,domain=-1:1,axis lines=middle,ticks=none, axis equal
    ,ymax=1.25,xmax=1.5,view={0}{90}
,xlabel={$u_1$},ylabel={$u_2$}]
\addplot+[no marks,thick,forget plot]{x} node {equilibria};
\ifnum#1>0
\addplot3+[-stealth,samples=6,no marks,thin
    ,quiver={u=-x+y,v=-y+x,scale arrows=0.2,every arrow/.append style={thin}}] {0};
\fi
\ifnum#1>1
\addplot+[no marks,thick] ({-x},{-x+x^3/3}) node[right] {centre manifold};
\fi
\end{axis}
\end{tikzpicture}
}}%

\item Neglect \(\partial_x\), \(u_1=u_2={}\)constant are equilibria---plotted in margin.
\smi0

\item Neglect \(\partial_x\), see \(\uv=(u_1,u_2)=e^{\lambda t}\vv\) leads to eigen-problem: \ldots \begin{itemize}
\item \(\lambda=0\) corresponding to \(\uv=(u_1,u_2)\propto(1,1)\), and \item \(\lambda=-2(a+b)\) corresponding to \(\uv\propto(1,-1)\).
\end{itemize}
So linearly, the subspace of equilibria \(\uv\propto(1,1)\) is attractive with rapid transients~\(e^{-2(a+b)t}\)---2nd plot in margin.
\smi1

\item Centre Manifold Theory asserts that under perturbation, such as here by small non-zero spatial gradients, the rapidly-attractive subspace is bent to a nearby rapidly-attractive manifold \footnote{A \emph{manifold} is a smooth curve, surface, \ldots, that we may parametrise with one, two, \ldots\ real parameters.  Read \protect\url{https://en.wikipedia.org/wiki/Center_manifold} for a flavour of the theory.}
on which there is slow evolution \cite[e.g.,][]{Carr81,Haragus2011}---3rd marginal plot.
\smi2

\item To construct the slow centre manifold, first \emph{choose} to parametrise it by the average \(U=(u_1+u_2)/2\) so that \(\uv\approx (1,1)U\).  Second, some algebraic machinations (\verb|homoDiff.txt|, \cref{AhomoDiff}) constructs the slow centre manifold to be
\begin{subequations}\label{eqssvhdsm2}%
\begin{align}&
u_{1,2}=U\pm\frac{a-b}{2(a+b)}\D xU+\Ord{\partial_x^3},
\label{eqsvhdsm2u}
\end{align}
on which the parameter~\(U\) evolves according to the \pde
\begin{equation}
\D tU=\frac{2ab}{a+b}\DD xU+\Ord{\partial_x^3}.
\label{eqsvhdpde2}
\end{equation}
\end{subequations}
I invite you to substitute to verify.
\end{enumerate}

\begin{activity}  For the specific case of diffusivities \(a=1\) and \(b=3\) substitute~\eqref{eqssvhdsm2} into the microscale \ode{}s~\eqref{eqssvhdmic} and verify the \ode{}s are satisfied to residuals~\Ord{\partial_x^3}.
\end{activity}

\paragraph{Question:} what does the time varying \(U(x)\) and its spatial derivatives mean when the underlying lattice is discrete in space?   
After all, such a function~\(U(x)\) does not contain any information about the \emph{phase} of the underlying lattice: the model \pde\ would be the precisely same if we shifted the lattice by any fraction of the lattice spacing. 
So the answer is that \cref{fig:dbl} does not go far enough in its doubling-up of the dynamical system.  
I contend the modelling only really makes sense when we regard it as modelling an \emph{infinite ensemble} of lattices, each lattice at a different phase shift in the lattice spacing.
Provided the ensemble is started in time in some consistent ensemble of initial conditions, then the field~\(U(x,t)\) will evolve smoothly and model the evolution of the whole ensemble.

\paragraph{Microscale information resolved}
The expression~\eqref{eqsvhdsm2u} informs us of the microscale structure.
Knowing this, and its higher-order refinements, empowers systematic modelling (cf.~Daniel Peterseim's lectures (DP), \S2.4).

\subsection{and underpins construction of higher order corrections}

Computer algebra routinely computes higher order models, even for complicated multiphysics scenarios.  
Here \verb|homoDiff.txt| constructs, in terms of the mean \(\am=(a+b)/2\) and difference \(\ad=(a-b)/2\),
\begin{align*}&
u_{1,2}=U\pm\frac{\ad}{2\am}\D xU 
\mp\left(\frac{\ad}{6\am}-\frac{\ad^3}{8\am^3}\right)\Dn x3U 
+\Ord{\partial_x^5},
\\&
\D tU=\left(\am-\frac{\ad^2}\am\right)\DD xU
+\left(\frac\am{12}+\frac{\ad^2}{6\am}-\frac{\ad^4}{4\am^3}\right)\Dn x4U 
\\&\qquad{}
+\Ord{\partial_x^5}.
\end{align*}

But we do not have to stop there.
For the example case of \(a=1\) and \(b=3\), I computed to 40th~order, then estimated \cite[]{Mercer90} that the power series for solutions \(U\propto\cos(kx)\) converges for wavenumbers \(k<2.5\)\,, equivalent to length scales \(\Delta x>1.2d\).
\begin{itemize}
\item That is, this homogenisation should be good for predicting structures down to just two to three lattice spacings!

\item Also, the modelling neglects transients \(e^{-2(a+b)t}=e^{-8t}\); that is, the model resolves time scales significantly longer than~\(1/8\), say \(\Delta t>1/2\). 

\end{itemize}
I expect similar \emph{quantitative} bounds for general diffusivities (most other approaches only engender \emph{qualitative} bounds).

\subsection{Boundary conditions}
\label{ssBCs}

A challenge is to derive correct boundary conditions for macroscale \pde\ models such as the diffusion~\eqref{eqsvhdpde2}.
A bigger challenge is what boundary conditions should be applied to higher-order versions.
Here we only answer the first, not the second: the general approach for both was developed some decades ago \cite[]{Roberts92c}, and only relatively recently applied to this sort of homogenisation \cite[]{Chen2014}.
For definiteness I analyse the case when diffusivities are \(a=1\) and \(b=3\), and leave the class of general period-two diffusivities to you.

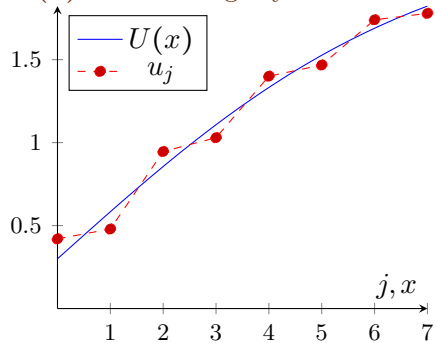
\begin{figure}
\centering
\caption{\label{fig:svhdbcs}given a boundary value at \(x=0\) of~\(u_0\) the microscale lattice values~\(u_j\) with diffusivities \(a=1\) and \(b=3\) zig-zags as shown.  Consequently the macroscale field~\(U(x)\) must be slightly different to~\(u_0\) at the boundary.}
\begin{tikzpicture}
\begin{axis}[small,domain=0:7,ymin=0,axis lines=middle,xlabel={$j,x$},legend pos=north west]
\addplot+[no marks]{0.3+2*tanh(x/7)};
\addplot+[dashed,samples=8,mark=*]{0.3+2*tanh(x/7)+0.12*cos(180*x)*exp(-x/7)};
\legend{$U(x)$,$u_j$}
\end{axis}
\end{tikzpicture}
\end{figure}%
\cref{fig:svhdbcs} illustrates the scenario near a boundary.
Suppose that the microscale boundary condition on the left, \(x=0\), is that~\(u_0\) is specified.
Then the microscale solution on the lattice must have a zig-zag structure, as shown, that reflects the alternating diffusivities, here starting on the left with the higher diffusivity \(b=3\).
Now we expect the macroscale \emph{mean} field~\(U(x)\) to travel smoothly through the `centre' of these zig-zags as illustrated.  
Consequently (\cref{fig:svhdbcs}), its boundary value~\(U|_{x=0}\) is generally different to the prescribed microscale value~\(u_0\).
We proceed to argue that the correct boundary condition for the diffusion \pde~\eqref{eqsvhdpde2} (\(a=1\) and \(b=3\)) is the Robin condition
\begin{equation}
U+\frac14\D xU=u_0\text{ at }x=0\,,
\label{eqsvhdbc2}
\end{equation}
dimensionally \(U+\frac d4\D xU=u_0\).
These express that in the scenario of \cref{fig:svhdbcs} \(U|_{x=0}\) must be a bit less than~\(u_0\).

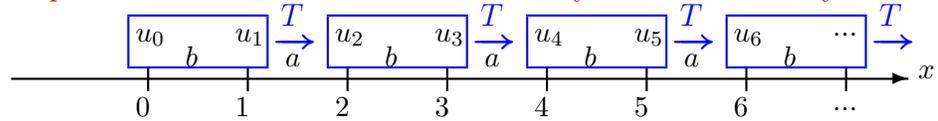
\begin{figure*}
\caption{schematic diagram of inhomogeneous diffusion on a microscale lattice of material~$u_j$: boundary conditions come from considering the map~\(T\) from one cell to the next further away from the boundary.}
\label{fig:svhldbc}
\centering
\setlength{\unitlength}{1ex}
\begin{picture}(75,9)
%\put(0,0){\framebox(75,9){}}
\put(0,-2){\thicklines
\put(0,5){\vector(1,0){72}$\ x$}
\multiput(11,4)(8,0){8}{\line(0,1)2}
\setcounter{i}{0}
\multiput(10,2)(8,0){7}{$\arabic{i}$\stepcounter{i}}
\put(66,2){$\cdots$}
\setcounter{i}{0}
\multiput(10,8)(8,0){7}{$u_{\arabic{i}}$\stepcounter{i}}
\put(66,8){$\cdots$}
\multiput(22,6)(16,0)3{$a$}
\multiput(14,6)(16,0)4{$b$}
\color{blue}
\multiput(9.5,6)(16,0){4}{\framebox(11,4){}\put(0.5,1){\LARGE$\xrightarrow T$}}
}
\end{picture}
\end{figure*}%
The key is to consider the \emph{spatial evolution} away from the boundary into the interior \cite[]{Roberts92c, Chen2014}.
Specifically, because the microscale has period-two, we consider the `dynamics' of the map from one pair of consecutive lattice points to the next pair: the map~\(T\) shown in \cref{fig:svhldbc}.
Many will recognise this~\(T\) as the map from one `cell' of the problem to the next `cell'.

In the slow manifold model the evolution is slow, hence time derivatives are small.
Consequently, to a useful approximation we neglect time variations in the development of boundary conditions. 
Such time variations could be incorporated, but they obfuscate the main issues, have only a small effect, so are neglected.

To find the map~\(T\) from one cell to the next, without loss of generality we just derive the map from~\((u_0,u_1)\) to~\((u_2,u_3)\) as all the others are the same.
For the specific diffusivities \(a=1\) and \(b=3\), and neglecting time derivatives, the original governing \ode{}s~\eqref{eqsvhdode} are
\begin{equation*}
0=\begin{cases}
1[u_{2}-u_1]+3[u_{0}-u_1]&\text{for }j=1,\\
3[u_{3}-u_2]+1[u_{1}-u_2]&\text{for }j=2.
\end{cases}
\end{equation*}
\begin{itemize}
\item Rearrange the first to give~\(u_2\) as a function of~\((u_0,u_1)\):
\(u_2=-3u_0+4u_1\)\,.
\item Rearrange the second to give~\(u_3\) via a function of~\((u_1,u_2)\):
\(u_3=-\tfrac13u_1+\tfrac43u_2
=-\tfrac13u_1+\tfrac43(-3u_0+4u_1)
=-4u_0+5u_1\)\,.
\item Combining these two gives the map
\begin{equation*}
\begin{bmatrix} u_2\\u_3 \end{bmatrix}
=\underbrace{\begin{bmatrix} -3&4\\-4&5 \end{bmatrix}}_T\begin{bmatrix} u_0\\u_1 \end{bmatrix}.
\end{equation*}

\end{itemize}
The same relation hold for all cells/pairs, hence \(T^n\)~determines the spatial evolution away from the boundary into the interior (here an overall linear dependence upon~\(n\)).

Here we derive the boundary condition just from the first two cells/pairs.
Now, the macroscale field at \(x=1/2\) is the average of~\(u_0\) and~\(u_1\), that is, \(U|_{1/2}=\tfrac12u_0+\tfrac12u_1\).
Correspondingly, the macroscale field at \(x=5/2\) is the average of~\(u_2\) and~\(u_3\), that is, \(U|_{5/2}=\tfrac12u_2 +\tfrac12u_3\).
But we have expressions for~\(u_2,u_3\) in terms of~\(u_0\) and~\(u_1\), giving \(U|_{5/2}= -\tfrac72u_0 +\tfrac92u_1\)\,.
Also, near the boundary the macroscale solution is very nearly linear, so to a good approximation the macroscale \(U|_{1/2}\approx U+\tfrac12U_x\) where the right-hand side is evaluated at \(x=0\), and similarly \(U|_{5/2}\approx U+\tfrac52U_x\)\,.
Consequently we form the two linear equations (the left-hand sides are evaluated at \(x=0\))
\begin{align*}&
U|_{1/2}:\quad U+\tfrac12U_x =\tfrac12u_0+\tfrac12u_1\,,
\\&
U|_{5/2}:\quad U+\tfrac52U_x =-\tfrac72u_0 +\tfrac92u_1\,.
\end{align*}
Eliminate the unknown microscale~\(u_1\) by subtracting the second from nine times the first to give that
\begin{equation*}
U+\tfrac14U_x=u_0 \text{ at }x=0\,,
\end{equation*}
and hence establish the claimed boundary condition~\eqref{eqsvhdbc2}.

The corresponding argument, but backwards in space, provides correct boundary conditions on a right-end boundary.

\begin{activity} 
What is the boundary condition for general~\(b\) (keep \(a=1\))?
\end{activity}

\paragraph{Comments}
\begin{itemize}
\item Curvature in the macroscale solution here only arises from an evolving out-of-equilibrium solution, so neglecting curvature in~\(U(x)\) is a consistent equivalent approximation as neglecting time derivatives.

\begin{SCfigure*}
\centering
\caption{\label{fig:resultBC} in diffusion through two coupled heterogenous strands, the microscale (crosses) exhibits boundary layers at each end \protect\cite[from][Fig.~3]{Chen2014}.  Classic  arguments give \textsc{bc}s for macroscale \pde{}s that incorrectly predict the red line.  Our approach caters for the microscale boundary layers and heterogeneity to give \textsc{bc}s that correctly predict the blue line.}
\includegraphics[width=0.7\textwidth, 
height=\textheight,  keepaspectratio]{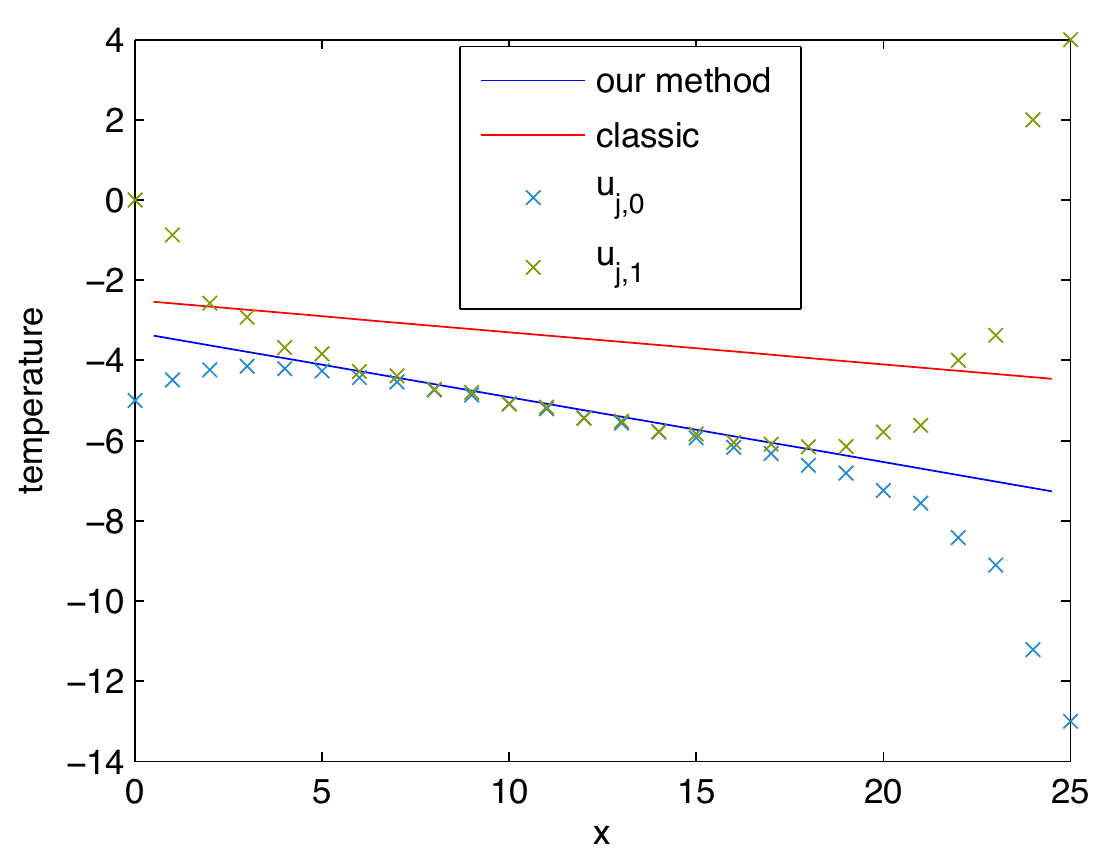}
\end{SCfigure*}%
\item In more complicated problems there are boundary layers at each boundary, shown for example by \cref{fig:resultBC}, that we cater for in a more general argument \cite[]{Roberts92c, Chen2014}.
\end{itemize}

\subsection{Microscale heterogeneous mechanical media}

What if the system illustrated by \cref{fig:svhld} is that of a mechanical system.
For simplicity, say there are unit masses at every microscale lattice point, connected by springs to its two neighbours.
The heterogeneity is that the springs are of alternating strengths~\(a\) and~\(b\).
We derive that the appropriate macroscale homogenisation is the classic wave \pde---with higher-order modification.

Let \(u_j(t)\) denote the displacement of each mass, and \(v_j(t)\)~denote the velocity of each mass.
Then the microscale \ode{}s are modified from~\eqref{eqsvhdode} by changing~``\(\dot u_j={}\)'' to ``\(\dot u_j=v_j\) and \(\dot v_j={}\)''.
As in \cref{ssrmrrs,ssepsurh}, non-dimensionalise and embed in an ensemble of phase shifts.
\begin{enumerate}
\item As before, focus on the slowly-varying solutions by regarding~\(\partial_x\) as `small'.
\item The eigenvalue equation become \(\lambda^2=0\) and \(\lambda^2=-2(a+b)\) indicating two slow modes among fast oscillations of frequency~\(\sqrt{2(a+b)}\).
Thus linearly, the slow subspace is \(\uv,\vv\propto(1,1)\) which \emph{acts as the centre of oscillations}, instead of being exponentially quickly attractive.
\item Theory asserts that under perturbation there exists an `asymptotically close' system that has a slow manifold, free of the fast oscillations, and tangent to the slow subspace \cite[\S2.5, this is a backward theorem!]{Roberts2018a}.
\item Almost the same algebraic machinations (\verb|homoVibr.txt|, \cref{AhomoVibr}) constructs the slow manifold to be
\begin{subequations}\label{eqssvhvsm2}%
\begin{align}&
u_{1,2}=U\pm\frac{a-b}{2(a+b)}\D xU+\Ord{\partial_x^3},
\\&
v_{1,2}=V\pm\frac{a-b}{2(a+b)}\D xV+\Ord{\partial_x^3},
\end{align}
on which the parameters~\(U,V\) evolves according to the \pde
\begin{equation}
\D tU=V,\quad \D tV=\frac{2ab}{a+b}\DD xU+\Ord{\partial_x^3}.
\label{eqsvhvpde2}
\end{equation}
\end{subequations}
I invite you to substitute to verify.
\end{enumerate}

The macroscale \pde~\eqref{eqsvhvpde2} is the classic wave \pde\ for the mechanical medium.
Straightforward higher-order analysis constructs higher-order models that show the waves are at least a little dispersive.

Boundary conditions for the wave \pde\ may be derived as in \cref{ssBCs}.

\paragraph{Fast waves may resonate}
However, in \emph{nonlinear} wave systems, one important difference is that the slow macroscale evolution is different when there are fast waves present, compared to when the fast waves are absent \cite[Ch.~13]{Roberts2014a}: the difference is typically quadratic in the fast wave amplitude (e.g., Stokes drift in water waves).
For example, consider the toy nonlinear system 
\begin{equation*}
\dot x=-x^3+(y^2+z^2)x,\quad \dot y=-\omega z,\quad \dot z=+\omega y\,.
\end{equation*}
Fast waves of frequency~\(\omega\) are in \((y,z)\), so the slow manifold is \(y=z=0\) exactly.  On this \(\dot x=-x^3\) always decays to zero.  But in the presence of fast waves the long-term solution is fundamentally different: put~\(y,z\) in polar coordinates~\(r,\theta\) and then \(\dot r=0\), \(\dot\theta=\omega\) and \(\dot x=(r^2-x^2)x\). 
So solutions with fast waves have \(x(t)\to\pm r\neq 0\) of the slow manifold prediction.
This mean effect of fast waves is independent of the frequency!

\begin{activity}
See a similar effect in \(\dot x=z^2\), \(\dot y=-z\) and \(\dot z=y\) via the coordinate transform that \(X=x+yz/2\), \(Y=y\) and \(Z=z\).
What is the overall \(x\)-evolution?
\end{activity}

\subsection{Optional: Nonlinear pattern formation is analogously rigorously supported}
\label{ssnpfars}

Recall that \cref{sammt} discussed that modelling an ensemble of phase shifted diffusivity was a rational way to form macroscale models of heterogeneous material.
A similar approach works when the heterogeneity is an emergent phenomena of the system.

Pattern formation is a common phenomena in science and engineering:
for examples, the stripes on a zebra, the spots on a leopard, and the ordered arrays of clouds.
Let's overview briefly one of the basic toy problems in this class.
Consider the small amplitude solutions of the Swift--Hohenberg system in one space dimension: a field~\(u(x,t)\) satisfies the nondimensional nonlinear `microscale' \pde
\begin{equation}
\D tu=ru-(1+\partial_{xx})^2u-u^3
\label{eq:shpde}
\end{equation}
on a domain of large extent in~\(x\).
For parameter \(r\)~small, the slow marginal modes are \(u\propto e^{\pm ix}\).
The aim is to derive, as a macroscale model over large~\(x\), the well-known Ginzburg--Landau \pde
\begin{equation}
\D tc\approx rc-3|c|^2c+4\DD x c\,,
\label{eq:gle}
\end{equation}
governing the complex amplitude~\(c(x,t)\) of oscillatory patterns \(u(x,t)\approx ce^{ix}+\bar ce^{-ix}\) \cite[e.g.]{Cross93}.

\begin{figure*}
\caption{cylindrical domain of the embedding \pde~\eqref{eq:shembed} for field \(\fu(\fx,y,t)\).  Obtain solutions of the Swift--Hohenberg \pde~\eqref{eq:shpde} on the blue line as \(u(x,t)=\fu(x,x+\phi,t)\) for any constant phase~\(\phi\) \protect\cite[Fig.~1]{Roberts2013a}.}
\label{fig:pattensemble}
\centering
\setlength{\unitlength}{0.01\linewidth}
\begin{picture}(96,27)
%\put(0,0){\framebox(96,26){}}
\put(0,5){\vector(1,0){96}}
\put(40,1){\(x\)}
\put(5,0){\vector(0,1){26}}
\put(2,25){\(y\)}
\put(0,20){\(2\pi\)}\put(4,21){\line(1,0){2}}
\put(2,5){\(0\)}
\put(30,22){domain \(\XX\times[0,2\pi)\)}
\put(30,12){\(\fu (\fx,y,t)\)}
\thicklines
\put(5,5){\framebox(88,16){}}
\color{blue}
\multiput(13,5)(16,0)5{
  \put(0,0){\line(1,1){16}}
  \multiput(0,0)(0,2)8{\line(0,1){1}}
  }
  \put(13,21){\line(-1,-1){8}}
  \put(2,12){\(\phi\)}
\color{magenta}
\put(70,0){\line(0,1){25}\,\color{black}\(\fx=X\)}
\put(70,13){\oval(10,20)}
\end{picture}
\end{figure*}
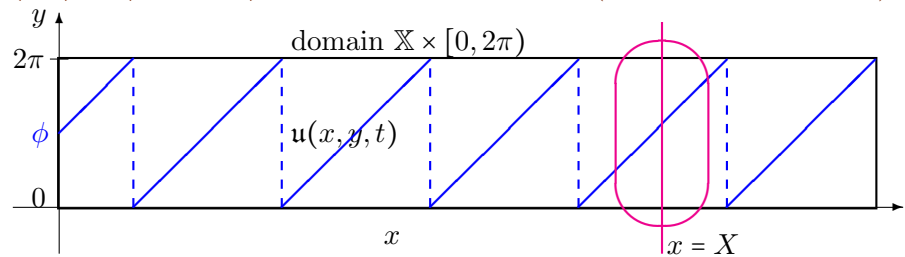
The trick to a rigorous approach is to embed the \pde~\eqref{eq:shpde} in the larger problem of analysing the ensemble of all phase shifts of the pattern \cite[\S2.5, \S3.3]{Roberts2013a}.
As indicated schematically in \cref{fig:pattensemble},
and in terms of a new ensemble\slash phase variable~\(y\), consider a new field~\(\fu (\fx,y,t)\) satisfying the \pde
\begin{equation}
\D t\fu =r\fu-(1+\partial_{yy}+2\partial_{y\fx}+\partial_{\fx\fx})^2\fu -\fu^3,
\label{eq:shembed}
\end{equation}
for \((\fx,y)\in \XX\times[0,2\pi)\),
where the field~\(\fu \) is \(2\pi\)-periodic in~\(y\).
Given any solution~\(\fu \) of the \pde~\eqref{eq:shembed}, elementary calculus shows that, for any chosen fixed phase~\(\phi\) and using that \(\fu \)~is \(2\pi\)-periodic in~\(y\), the field \(u(x,t)=\fu (x,x+\phi,t)\) (along the blue lines in \cref{fig:pattensemble}) is a solution of the Swift--Hohenberg \pde~\eqref{eq:shpde}.
Thus modelling of the dynamics of the ensemble \pde~\eqref{eq:shembed} immediately leads to models for the dynamics of the Swift--Hohenberg \pde~\eqref{eq:shpde}.

\emph{This `embedding' approach immediately makes new sense of the multiple space scales that others introduce, such as Elliot Carr's \xv\ and~\(\yv/\epsilon\).}  The `cell problem' here becomes solving in for the \(y\)~structure, given slow variations in~\(x\), but now with a well defined geometry.

The macroscale modelling of the ensemble \pde~\eqref{eq:shembed} may be done via rigorous \emph{local} Taylor expansions about an arbitrary station~\(x=X\) (\cref{fig:pattensemble}).  Such local models are typically only weakly coupled to neighbouring locales, and so the collection of local models generates a global \pde\ as the macroscale model  \cite[]{Roberts2013a}. 
Here we would derive the Ginzburg--Landau \pde~\eqref{eq:gle}.

The rigorous ensemble embedding here replaces heuristic multiple space and time scale assumptions traditionally employed in the asymptotic analysis of patterns \cite[e.g.]{Cross93, Vandyke87}.

\subsection{Open problems}

\begin{itemize}
\item Provide tools to automatically construct such models and boundary conditions for users, tools analogous to the web services I operate \cite[e.g.]{Roberts07d}.

\item Deduce quantitative bounds on the spatio-temporal resolution of many of the classic macroscale \pde\ closures in multiscale systems (mainly only for linear problems).

\item What are appropriate initial conditions \cite[e.g.]{Roberts89b} for a macroscale \pde\ given that the microscale is spatially discrete, and we only derive the macroscale by an embedding ensemble?

\item Develop Backwards Theory for centre\slash stable\slash unstable\slash slow\slash fast manifolds: that is, instead of ``for a given system provided restrictions there exists \ldots'' establish ``generally there exists a nearby system for which exact manifolds are \ldots''
\cite[aka][]{Grcar2011} \cite[]{Roberts2018a}.

\item Rigorously support and practical procedures to develop boundary conditions for macroscale \pde{}s in 2D or 3D spatial domains.

\item What if the microscale is stochastic? \cite[may inspire an approach]{Roberts06k}
\end{itemize}

\section{Multiscale computation of microscale systems}
\label{smcms}
\localtableofcontents

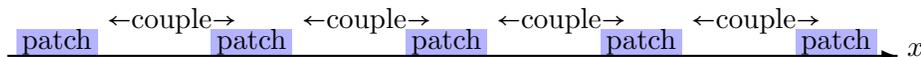
\begin{figure*}
\centering
\caption{\label{fig:coupat}schematic of spatial patches separated by unsimulated space.}
\setlength{\unitlength}{0.01\textwidth}
\begin{picture}(100,6)
%\put(0,0){\framebox(100,6){}}
\thicklines
\multiput(10.8,3)(21,0)4{\(\leftarrow\)couple\(\to\)}
\multiput(1,1)(21,0)5{\colorbox{blue!30!white}{\makebox(7,1){patch}}}
\put(0,0){\vector(1,0){96}$\ x$}
\end{picture}
\end{figure*}

Suppose that in some problem 
\begin{itemize}
\item we have a trustworthy microscale simulator, 
\item and a `spatial' domain so large the microscale code is not feasible,
\item but we do \emph{not} know and/or \emph{cannot} derive a macroscale closure.
\end{itemize}
Answer: use the microscale simulator on small patches of space (\cref{fig:coupat}), with the patches coupled over unsimulated space, craftily, so that we make macroscale predictions.
Such predictions are computed relatively quickly when the patches are a small fraction of the whole domain. 
\cref{sscicpc,sswdss} show that the macroscale homogenisation of the coupled patches is accurate in cases when we do know the macroscale closure.

\subsection{One Patch to rule them all, \ldots}

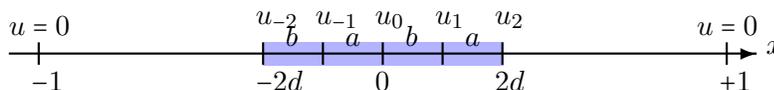
\begin{figure*}
\centering
\caption{\label{figoprta}one `small' patch of heterogeneous diffusion, in the centre, to be coupled to the `distant' boundary values of \(u=0\) at \(x=\pm1\).}
\setlength{\unitlength}{0.6ex}
\begin{picture}(103,12)
%\put(0,0){\framebox(103,12){}}
\thicklines
\put(34,4.5){\colorbox{blue!30!white}{\makebox(30,1){}}}
\put(0,5){\vector(1,0){100}$\ x$}
\put(15,0){
\multiput(19,3.5)(8,0){5}{\line(0,1)3}
\multiput(-11,3.5)(92,0){2}{\line(0,1)3\put(-4,4){$u=0$}}
\put(18,0){$-2d$} \put(34,0){$0$} \put(50,0){$2d$}
\setcounter{i}{-2}
\multiput(18,9)(8,0){5}{$u_{\arabic{i}}$\stepcounter{i}}
\multiput(22,6)(16,0)2{$b$}
\multiput(30,6)(16,0)2{$a$}
\put(-12,0){$-1$}
\put(80,0){$+1$}
}
\end{picture}
\end{figure*}%
The simplest scenario is just one patch coupled to distant boundary values.
\cref{figoprta} shows one small patch of length~\(4d\) in a domain \(-1<x<1\).
Inside the patch the microscale is that of heterogeneous, period-two, diffusion.
This microscale patch is to be coupled to the boundary values at the distant \(x=\pm 1\) of \(u=0\).
How can we couple to make correct predictions?

Here we know the desired predictions.
The macroscale closure is that \(U_t=DU_{xx}\) (approximately) for effective homogenised diffusivity \(D=2ab/(a+b)\).
So that, for example, if the macroscale field is approximately parabolic, \(U(x,t)=U_0(t)(1-x^2)\), then \ldots\ the macroscale closure gives \(\dot U_0=-4ab/(a+b)\,U_0\).  
We compare the patch scheme with this \ode.

Within the microscale patch, the microscale equations, from~\eqref{eqsvhdode}, in the interior of the patch are
\begin{subequations}\label{eqsoptra}%
\begin{align}&
d^2\dot u_{-1}= a[u_{0}-u_{-1}]+b[u_{-2}-u_{-1}],
\\&
d^2\dot u_{0}= b[u_{1}-u_0]+a[u_{-1}-u_0],
\\&
d^2\dot u_{1}= a[u_{2}-u_{1}]+b[u_{0}-u_{1}],
\end{align}
where the \(d^2\)-factor on the left-hand side caters for the microscale lattice spacing of~\(d\). 
Whereas \cref{sammt} non-dimensionalises the lattice spacing to \(d=1\), here we scale the macroscale domain length to two and so the microscale length scale is denoted~\(d\).   

But we need two edge values for the patch, at \(x=\pm2 d\).
These come from macroscale coupling with the boundaries---although in general from coupling with neighbouring patches.
Let's pose the predicted macroscale is classic parabolic interpolation through the two given boundary values, namely zero at \(x=\pm1\), and the evolving centre-patch value of~\(u_0(t)\) at the centre-patch \(x=0\).
That is, the predicted macroscale is \(U=(1-x^2)u_0(t)\).  
Then the two patch-edge values are taken to be this field evaluated at the edges \(x=\pm 2d\), namely the microscale 
\begin{equation}
u_{\pm2}=(1-4d^2)u_0\,.
\label{eqoptraev}
\end{equation}
\end{subequations}
This coupling completes the patch simulation equations.

What are the dynamics of the patch scheme? 
It is linear so find a general solution, after substituting~\eqref{eqoptraev}, via seeking \((u_{-1},u_0,u_1)=e^{\lambda t}\vv\)\,.
That is, find the eigen-values and eigen-vectors of matrix
\begin{equation*}
\frac1{d^2}\begin{bmatrix} -(a+b)&a+b-4bd^2&0\\
a&-(a+b)&b\\ 0&a+b-4ad^2&-(a+b) \end{bmatrix}.
\end{equation*}
In terms of \(\mu=\lambda d^2\), this matrix's characteristic polynomial factors to
\begin{equation*}
(\mu+a+b)\left[\mu^2+2(a+b)\mu+8abd^2\right]=0\,.
\end{equation*}
That is, \(\mu=-(a+b)\) and \(\mu=-(a+b)\pm\sqrt{(a+b)^2-8abd^2}\).
In the interesting cases of small patches, small~\(d\), these give the three eigenvalues
\begin{equation*}
\lambda=\frac\mu{d^2}\approx \frac{-4ab}{a+b},\ -\frac{a+b}{d^2},\ -\frac{2(a+b)}{d^2}.
\end{equation*}
For small~\(d\) the last two are large negative eigenvalues and the first is relatively small.
Thus, quickly, after a cross-patch diffusion time of~\(d^2/(a+b)\), all patch solutions lie on the subspace corresponding to the small eigenvalue.
On this slow subspace they evolve nearly as~\(\exp[-4ab/(a+b)\,t]\).
This rate,~\(-4ab/(a+b)\), exactly matches the macroscale closure \pde.
The patch scheme gives good macroscale predictions despite solving the microscale on only the small fraction~\(2d\) \text{of the domain}.

\begin{activity} What does the patch scheme predict when the boundary condition at \(x=-1\) is replaced by \(U_x=0\)?  
For simplicity set diffusivities \(a=1\) and \(b=3\)\,.  
How does the prediction compare with the macroscale homogenisation?
\end{activity}

\subsection{Automatic macroscale closure}

I think this amazing!
We only solved the microscale equations on a patch of length~\(4d\) and yet predict the macroscale homogenisation.
Similarly for executing microscale code instead of solving.
The patch size~\(4d\) could be tiny and we still get an accurate homogenisation.
We have not coded into the patch scheme any knowledge of the homogenisation, yet the scheme effectively, and on-the-fly, forms an accurate macroscale closure. 

There are two caveats here.
Firstly, in the scheme analysed we only predict one macroscale mode over the domain, the parabola~\(1-x^2\).
But what about other macroscale modes such as \(\sin(\pi x)\) and \(\cos(3\pi x/2)\)?
To resolve all three of these modes we would need three patches in the domain.
In general, \(N\)~patches resolve \(N\)~macroscale modes (\cref{sscicpc}).

\begin{table}
\caption{\label{tbloph}effective diffusivity for the macroscale mode with one patch on the microscale of~\([-nd,nd]\) when the microscale has period-two: even~\(n\) is always accurate.}
\begin{equation*}
\begin{array}{rl}
\hline
n&\text{effective diffusivity}
\\\hline
1&(a+b)/2
\\2&2ab/(a+b)
\\3&\frac92ab(a+b)/(2a^2+5ab+2b^2)
\\4&2ab/(a+b)
\\5&\frac{25}2ab(a+b)/(6a^2+13ab+6b^2)
\\6&2ab/(a+b)
\\7&\frac{49}2ab(a+b)/(12a^2+25ab+12b^2)
\\8&2ab/(a+b)
\\\text{odd}&2ab(a+b)/\big[(a+b)^2-(a-b)^2/n^2\big]
\\\hline
\end{array}
\end{equation*}
\end{table}%
Secondly, let's check what happens with different number of lattice points in the patch.
Computer algebra, \verb|onePatchHomo.txt| (\cref{AonePatchHomo}), easily checks for patches of size~\([-nd,nd]\).
All cases have \begin{itemize}
\item \((2n-2)\) negative eigenvalues of large magnitude\({}\propto 1/d^2\) characteristic of rapid sub-patch diffusion, and 
\item one small negative eigenvalue giving the effective macroscale homogenised diffusivities of \cref{tbloph}.
\end{itemize}

Patches of size with even~\(n\) are macroscale correct (and odd~\(n\) are increasingly accurate as \(n\)~increases).
This illustrates the general rule that it is best if the half-size of a patch is an integral multiple of the microscale periodicity \cite[DP, and][\S5.2]{Bunder2013b}.

\begin{activity}
Consider DP's first problem but over \(-1<x<1\): \(-[A_\epsilon(x)u_x]_x=1\) with \(u(\pm 1)=0\) and \(A_\epsilon=[2+\cos(2\pi x/\epsilon)]^{-1}\) for integer~\(1/\epsilon\).
Form one patch \(-\epsilon<x<\epsilon\) with edge conditions from the macroscale parabolic interpolation \(u(\pm\epsilon)=u(0)(1-\epsilon^2)\).
Solve exactly within the small patch (aided by symmetry) to discover the solution in the patch is that of the exact whole domain solution!
\end{activity}

\subsection{Optional: Nonlinear diffusion in one patch}
\label{ssndop}

How does a patch of nonlinear dynamics perform?
Let's consider the nonlinear \pde\ example for field~\(u(x,t)\) of 
\begin{subequations}\label{eqsndop}%
\begin{equation}
\D tu=u\DD xu\,,
\label{eqndop}
\end{equation}
subject to boundary conditions \(u=0\text{ at }x=\pm 1\)\,.

Solve only on the small patch \(|x|\leq h\).
The macroscale interpolated field is, given the centre-patch value \(U_0(t):=u(0,t)\)  that \(U(x,t)=(1-x^2)U_0(t)\).
So the edge values on the patch are that
\begin{equation}
u(\pm h,t)=(1-h^2)U_0=(1-h^2)u(0,t).
\label{eqndpev}
\end{equation}
\end{subequations}
What are the predictions when we only solve (compute) \pde~\eqref{eqndop} on the small patch?

\begin{enumerate}
\item \emph{Embed} the problem in a family of problems parametrised by~\(\gamma\), \(0\leq\gamma\leq1\), through generalising the edge condition~\eqref{eqndpev}:
\begin{equation}
u_t=uu_{xx},\quad u(\pm h,t)=(1-\gamma h^2)u(0,t).
\label{eqndopr}
\end{equation}
We use a theory based at \(\gamma=0\) to access results that hold for full coupling, \(\gamma=1\).
\pgfplotsset{ colormap={allred}{rgb=(0.7,0,0) rgb=(1,0,0)}}
\def\smi#1{\marginpar{
\begin{tikzpicture}
\begin{axis}[footnotesize,font=\footnotesize
    ,xmin=0,ymin=0,zmin=-1,xmax=1.1,ymax=1.2,zmax=1.1, ticks=none
    ,axis lines=middle, view={110}{25}
    ,xlabel={$\gamma$},ylabel={$U_0,V$},zlabel={$u(x)-U_0,W$}]
\addplot3+[no marks,very thick,forget plot,domain=0:1.1]({0},{x},{0}) 
    node[below left] {equilibria};
\ifnum#1>0
\addplot3+[-stealth,samples=4,no marks,thin,domain=0:1.1,y domain=-1:1
    ,quiver={u=0,v=0,w=-z,scale arrows=0.5,every arrow/.append style={thin}}] 
    ({0},{x},{y});
\fi
\ifnum#1>1
\addplot3+[surf,no marks,domain=0:1,y domain=0:1.1,opacity=0.5
    ,samples=7,colormap name=allred] {y*x^2/2};
\fi
\end{axis}
\end{tikzpicture}
}}%

\item \emph{Equilibria?}  \(u={}\)constant and \(\gamma=0\)\,.  Call the constant~\(U_0\) for compatibility---see marginal plot.
\smi0

\item \emph{Linearise:}  seek \(u=U_0+\hat u(x,t)\) for small~\(\hat u\) and negligible~\(\gamma\). 
Then the \pde\ and edge conditions become \ldots
\begin{equation*}
\hat u_t=U_0\hat u_{xx},\quad \hat u(\pm h,t)=\hat u(0,t).
\end{equation*}
Seek solutions via separation of variables, \(\hat u=e^{\lambda t}v(x)\), and find \ldots\  eigenvalues and eigenfunctions are
\begin{align*}&
\lambda_k=-U_0\pi^2k^2/h^2,
\\&
v_k(\xi)=\begin{cases}
\sin(k\pi x/h)&k=1,2,3,\ldots\,,\\
\cos(k\pi x/h)&k=0,2,4,\ldots\,,
\end{cases}
\end{align*}
as well as some generalised eigenfunctions for \(k=2,4,6,\ldots\)---schematic 2nd marginal plot.
\smi1

\item \emph{Emergence theory?} Because there is a zero eigenvalue with the rest negative and\({}\leq-U_0\pi^2/h^2\), centre manifold theory \cite[e.g.]{Carr81, Haragus2011} asserts the system~\eqref{eqndopr} possesses a 2D slow centre manifold~\(u(x,t)=u(x,U_0,\gamma)\) on which the system evolves, \(U_{0t}=g(U_0,\gamma)\). 
Further, the solutions on the slow centre manifold attract \emph{all} nearby ones roughly as~\(\exp(-U_0\pi^2t/h^2)\)---schematic 3rd marginal plot with red centre manifold.
\smi2

Consequently, provided evaluation at \(\gamma=1\) is valid, the original patch system~\eqref{eqsndop} possesses a quickly attractive 1D slow centre manifold~\(u=u(x,U_0,1)\) on which the system evolves \(U_{0t}=g(U_0,1)\). 
Again, \emph{on a cross-patch diffusion time all solutions of the patch system~\eqref{eqsndop} approach a 1D \ode\ which turns out to be the appropriate macroscale dynamics.}

\item \emph{Construction:} I chose this nonlinear problem as it is straightforward to verify the exact slow centre manifold is precisely
\begin{equation*}
u=(1-\gamma x^2)U_0\,,
\quad\text{s.t. }
U_{0t}=-2\gamma U_0^2\,.
\end{equation*}
To verify, substitute into~\eqref{eqndopr} \ldots.%
\footnote{This example also nicely illustrates two ways centre manifold models `break down'.  For parameter \(\gamma<0\) the centre manifold is attractive, but solutions within the centre manifold explode to infinity in finite time, via \(\dot U_0=2(-\gamma)U_0^2\).
Whereas for parameter \(\gamma>1\) the centre manifold exists and solutions within it are stable for all time, but the centre manifold surely no longer attracts all nearby solutions as the nonlinear diffusion coefficient,~\(u\), is negative for \(|x|>1/\sqrt\gamma\).}

Hence, at \(\gamma=1\), the original patch system~\eqref{eqsndop} for~\(u(x,t)\) has the attractive exact slow centre manifold\begin{equation}
u=(1-x^2)U_0\,,
\quad\text{s.t. }
\dot U_0=-2U_0^2\,.
\label{eqsndopsm}
\end{equation}
\end{enumerate}

The patch slow centre manifold~\eqref{eqsndopsm} also happens to be exact for the whole domain \pde~\eqref{eqndop}.
But here we discover the exact dynamics, \(\dot U_0=-2U_0^2\), `economically' by solving the \pde\ only on a microscale patch, and using interpolation over unsolved space to fill in the macroscale gaps.

\begin{activity}
What (little) would change in this analysis if the macroscale boundary conditions are \(2bu\pm(1-b)u_x=0\) at \(x=\pm1\)?  where parameter \(0\leq b\leq1\).
\end{activity}

\subsection{A basic atomic simulation}

So far I have discussed several toy problems.
Let's have a brief look at the realistically complicated scenario of atomistic simulation, and a patch simulation.

Consider a long domain \(|x|<H\), with thin square cross-section \(|y|,|z|<h\), filled with a monatomic gas.
We want to simulate the macroscale diffusion of heat along the long thin domain.
Here we know the macroscale is the diffusion \pde\ \(\D tT=D\DD xT\), although maybe not know the diffusivity~\(D\), nor its temperature\slash pressure dependence. 
But let's pretend we do not know even the form of this macroscale closure.
Instead, let's make a macroscale prediction for the diffusion of heat using a microscale atomic simulation within an single patch \cite[]{Alotaibi2017a}.

\begin{figure*}
\centering
\caption{trajectories of \(64\) atoms, over a time \(0\leq t\leq3\)\,, in a triply-periodic, cubic, spatial domain, showing the beginnings of the complicated inter-atomic interactions \protect\cite[Fig.~1]{Alotaibi2017a}.
View this stereo pair cross-eyed for a 3D effect.}
\label{fig:sim64}
\includegraphics[scale=0.85]{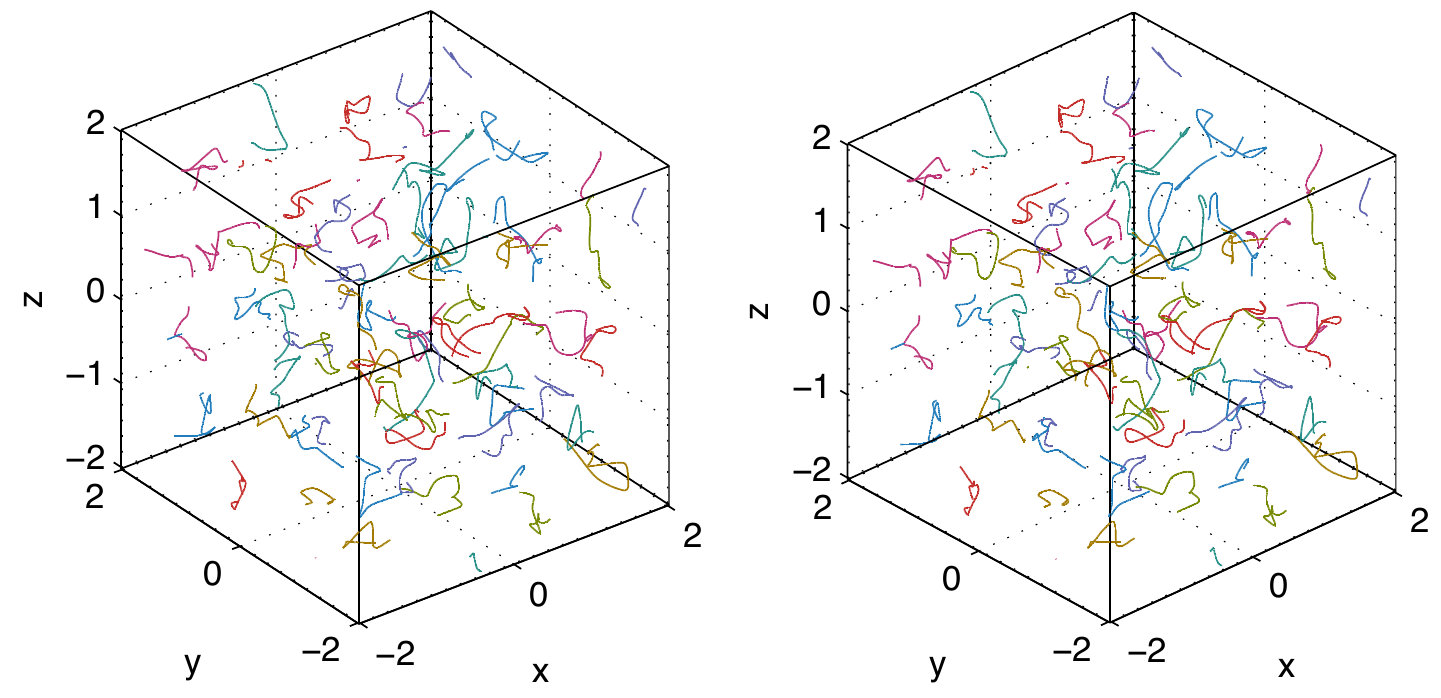}
\end{figure*}%
The \emph{most straightforward} atomic simulation to code is that of the motion of atoms, with interatomic forces determined from the classic Lennard-Jones potential, in a triply periodic, cubic, domain.
\cref{fig:sim64} shows such a simulation with \(64\)~atoms over a short microscale time.

Similarly,  Lattice--Boltzmann simulations easiest to code when periodic (as commented by CSF).

\begin{figure*}
\centering
\caption{the simplest case is one triply-periodic patch of atomistic simulation, \(-h<x<h\)\,, coupled to distant sidewalls, at \(x=\pm H\)\,, of specified temperature.  
The patch's core region defines its local temperature, and a proportional controller applied in the left and right action regions engenders a good macroscale prediction \protect\cite[Fig.~3]{Alotaibi2017a}.}
\label{fig:schemePPatch}
\includegraphics[scale=0.85]{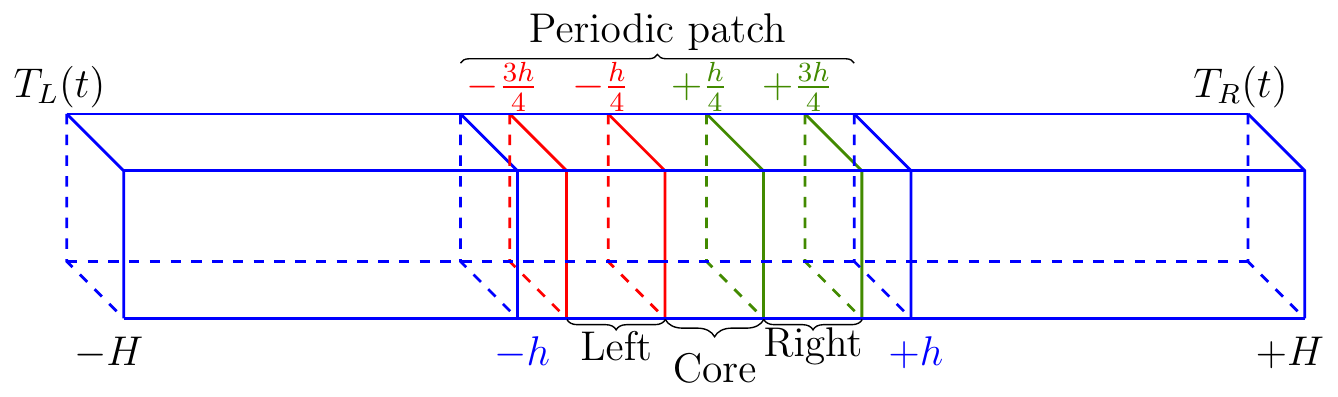}
\end{figure*}%
Question: how can we use such a microsimulation code in a patch scheme?
\cref{fig:schemePPatch} shows the simplest patch scheme.
We position just one micro-domain, a cube of side~\(2h\),  centred at \(x=0\), and filling the square cross-section of the domain.
The distant boundaries of the macro-domain are at \(x=\pm H\), on which we impose temperatures~\(T_L,T_R\).

But the micro-code patch is triply periodic, so we cannot specify boundary values on the edge of the patch, because the patch has no edge!
Instead we \emph{control} the patch.
Define four equi-sized regions in the patch as shown in \cref{fig:schemePPatch}:
\begin{itemize}
\item the core region---we estimate the macroscale temperature~\(T_0\) at \(x=0\) from the kinetic energy of the atoms in the core;
\item two action regions in which we apply a proportional controller to heat\slash cool the atoms depending upon whether the macroscale interpolated temperature through~\(T_L,T_0,T_R\) is less\slash more than the kinetic energy of the atoms in the particular action region;
\item and an `unmentioned' region whose role is to complete the microscale periodicity.
\end{itemize}

\begin{figure*}
\centering
\caption{(a)~temperatures over macroscale times in the sub-patch regions.
Simulate \(343\)~atoms in a patch of spatial periodicity~\(2h=7\) and with control strength \(\mu=30\) to couple with macroscale boundary temperatures \(T_R=1.5\) and \(T_L=0.5\) at \(x=\pm7\) \protect\cite[Fig.~5(a)]{Alotaibi2017a}.}
\label{fig:controlHeat}
\begin{tabular}{c@{}c}
\rotatebox{90}{\hspace{3ex}temps \(T_l\), \(T_0\), \(T_r\)} &
\includegraphics[scale=0.85]{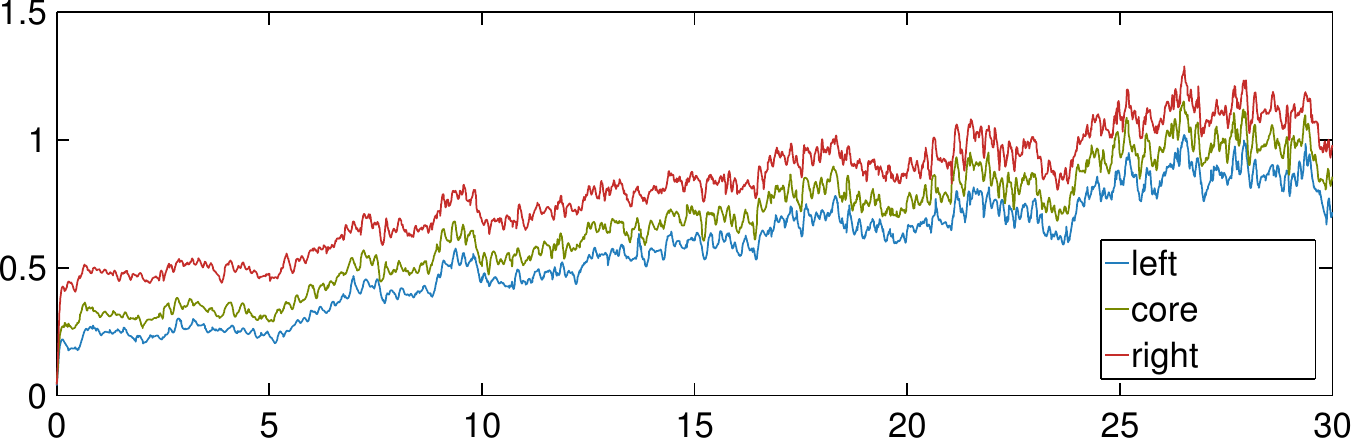}\\
(a)& time $t$
\end{tabular}
\end{figure*}%
\cref{fig:controlHeat} shows the result.
There is an initial equilibration transient which occurs on a time-scale almost too small to see on this time-axis.
After that rapid transient, a relatively cool patch, temperature about~\(0.3\), gradually heats up, with more heat flux from the right-end at \(T_R=1.5\), than from the left-end at \(T_L=0.5\)\,.
Over the macroscale time shown, the atoms in the patch heat up to the equilibrium temperature \(T_0\approx1\), with fluctuations due to the microscale chaos in the atomic motion.
Because the macroscale equilibrium temperature should be linear from left to right, in the final atomic `equilibrium' the action regions have temperatures \(T_L<T_0<T_R\), by roughly \text{equal amounts.} 

\needspace{4\baselineskip}
This controlled periodic-patch scheme does appear to predict reasonably correctly the macroscale dynamics.
We analysed the dynamics of such a single coupled periodic patch, and multiple coupled periodic patches, to determine optimal control parameters \cite[\SS4,5]{Alotaibi2017a}.

\subsection{Classic interpolation couples patches consistently,} 
\label{sscicpc}

Let's turn now to the issue of coupling the computation on many microscale patches across a large macroscale domain.
The simple answer is to couple by providing edge values for the patch computation via classic Lagrangian interpolation over the macroscale gaps of the patch centre-values \cite[e.g.]{Roberts06d}.
We show that the homogenisation of the dynamics of patches and gaps is reasonably accurate. 

\begin{SCfigure*}
    \centering
	\caption{Gap-tooth solution of Burgers' \pde\
	on~$[0,2\pi]$ through microsimulation on eight patches, each of small
	width; the teeth are coupled by classic Lagrangian interpolation (from \cref{sspfeft}).}
    \label{fig:burg3}
\includegraphics[width=0.76\textwidth, 
height=\textheight,  keepaspectratio]{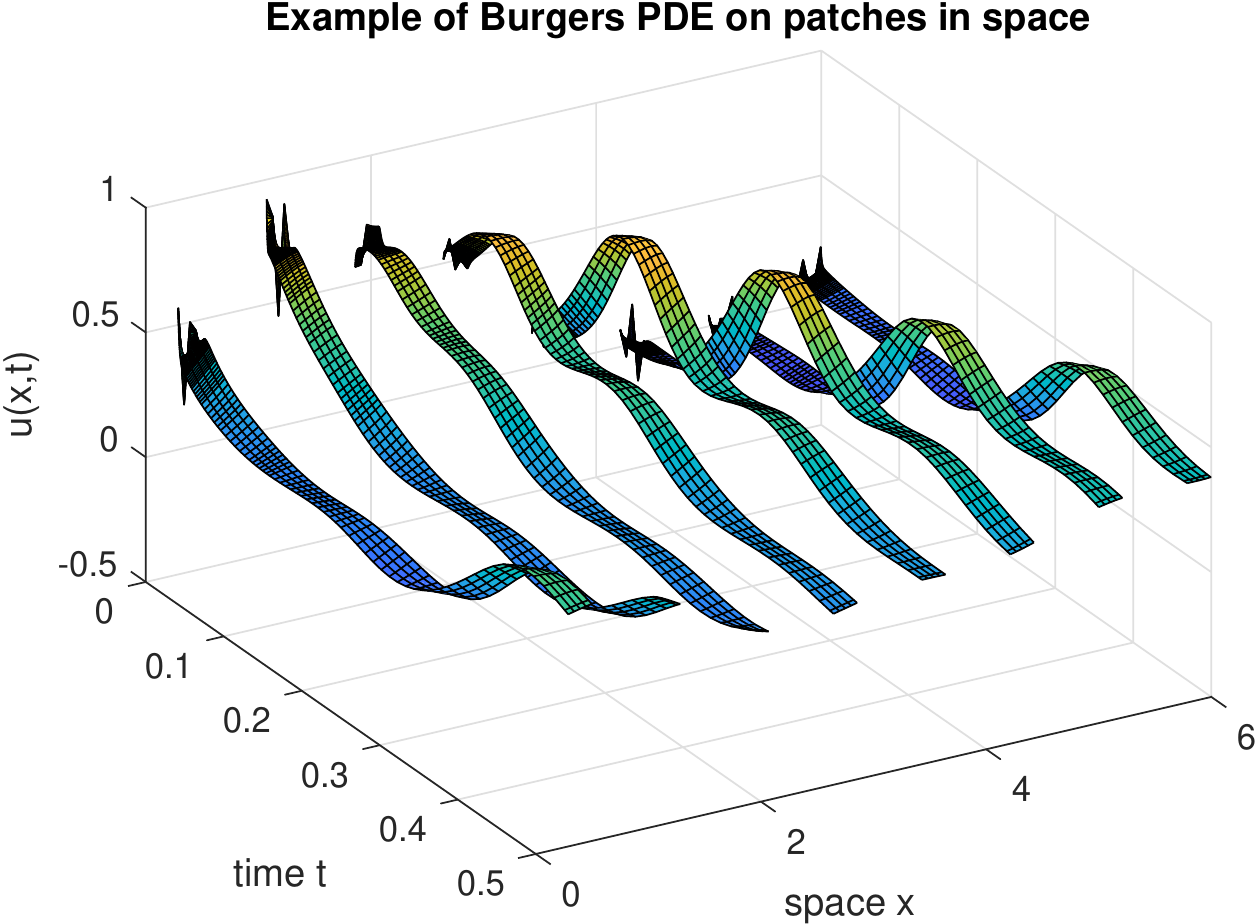}
\end{SCfigure*}%
\cref{fig:burg3} illustrates an example applied to the field~\(u(x,t)\) satisfying Burgers' \pde\ \(u_t+30uu_x=u_{xx}\).
We only compute on the patches, and not at all in the gaps.
Initial rapid transients, only just visible in the initial instants, decay to smooth sub-patch fields that then interact and evolve over macroscale space-times.

\paragraph{Equation-Free Toolbox} 
\footnote{\protect\url{https://github.com/uoa1184615/EquationFreeGit.git}}
Execute \verb|patchConfig1()| to see \cref{fig:burg3} generated by the example code near the start of the function.

In the macroscale domain place a grid, spacing~\(H\), with grid points~\(X_j\).
Centre a microscale patch, of size~\(2h\), at each grid point so the \(j\)th~patch covers \(X_j-h\leq x\leq X_j+h\).
Let~\(u_j\) denote the field in the \(j\)th~patch, but often it is convenient to use the \(j\)th~patch-centric space variable \(\xi=(x-X_j)/H\), \(-r\leq\xi\leq r\) for scale ratio \(r=h/H\), to describe sub-patch structures, so usually we consider~\(u_j(\xi,t)\).
The macroscale field is then formed by interpolating over the gaps all the centre-patch values \(U_j(t):=u_j(0,t)\).

To create the inter-patch coupling, define the macroscale shift operator $Eu(x):=u(x+H)$ and equivalently $EU_j:=U_{j+1}$, as appropriate for steps on the coarse grid size~$H$, and then its inverse gives \(E^{-1}u(x)=u(x-H)\) and $E^{-1}U_j:=U_{j-1}$.  Consequently, \(E^{\pm p}u(x)=u(x\pm pH)\) and \(E^{\pm p}U_j=U_{j\pm p}\) which naturally holds for all real~\(p\).
In particular, the patch edge-values need to be the interpolated macroscale: that is, since the patch-edges are at \(\xi=\pm r\) and patch-centre is at \(\xi=0\) for scale ratio \(r=h/H\)\,,
\begin{equation}
u_j(\pm r,t)=E^{\pm r}u_j(0,t)=E^{\pm r}U_j(t).
\label{eqcicpdev}
\end{equation}
The patch scheme is then to compute\slash solve the microscale code\slash\pde\ in each of the patches, \(|x-X_j|<h\), coupled by using~\eqref{eqcicpdev} to specify the patch edge-values where \(U_j\)~denotes the centre value of the \(j\)th~patch.

But how do we compute the required edge-values~\(E^{\pm r}U_j\)?
Answer: via some classic identities for discrete operators \cite[p.65, e.g.]{npl61}:
\begin{align*}
    \text{difference}\quad& \delta = E^{1/2}-E^{-1/2}\,,\\
    \text{mean}\quad& \mu = \half(E^{1/2}+E^{-1/2})=\sqrt{1+\rat14\delta^2}\,,\\
    \text{shift}\quad& E = 1+\mu\delta+\half\delta^2\,,
\end{align*}
Then, applied to~\(U_j\) for the edge values, and recall the scale ratio \(r=h/H\),
\begin{subequations}\label{eqsdoper}%
\begin{align}
    E^{\pm r}
    &=(1+\mu\delta+\half\delta^2)^{\pm r}
    \nonumber\\
    &=1 +r(\pm\mu\delta +\rat12r\delta^2 )
    \label{eqdopera}\\&\quad{}
    +r(r^2-1)(\pm \rat1{3!}\mu\delta^3
    +\rat1{4!} r\delta^4)
    \label{eqdoperb}\\&\quad{}
    +r(r^2-1)(r^2-4)(\pm\rat1{5!}\mu\delta^5
    +\rat1{6!}r\delta^6)
    \label{eqdoperc}\\&\quad{}
    +r(r^2-1)(r^2-4)(r^2-9)(\pm \rat1{7!}\mu\delta^7
    +\rat1{8!}r\delta^8)
    \label{eqdoperd}\\&\quad{}
    +\Ord{\delta^9}.
\end{align}
\end{subequations}
We approximate by truncating at some line: truncating at~\eqref{eqdopera} is locally second-order parabolic interpolation; truncating at~\eqref{eqdoperb} is locally fourth-order quartic interpolation; and so on.

A simple test of accuracy? Linear diffusion \(u_t=u_{xx}\), with macroscale periodicity of~\(2\pi\) and determine the accuracy of the macroscale modes.
The precise eigenvalues should be \(\lambda=-k^2\) for \(k=0,1,2,\ldots\) (except for \(k=0\), all are multiplicity two).
\cref{tblspatb} gives the numerical results and shows the macroscale eigenvalues, those for \(|k|<N/2\), are accurately determined by the patch scheme: the errors are \Ord{(kH)^4} as befits the fourth-order interpolation across gaps.
The rightmost column of \cref{tblspatb} gives the leading microscale eigenvalues which are several orders larger, corresponding to the decay of sub-patch modes on a sub-patch diffusion time~\(\Ord{1/h^2}=\Ord{N^2}\).
The patch scheme with multiple patches appears to make successful macroscale predictions, despite only computing on separated small patches of the domain.
\begin{table*}
    \centering
	\caption{eigenvalues~$\lambda$ of patch scheme modes for linear diffusion with $N$~patches, spaced $H=2\pi/N$,
	with scale ratio $r=0.1$\,,  $n=11$ points in each microscale patch; and with the fourth order coupling~\eqref{eqdoperb}  \protect\cite[from][]{Roberts06d}.}
    \label{tblspatb}
    \begin{tabular}{|r|llll|c|}
        \hline
        $N$ & \quad 1 & \quad 2,3 & \quad 4,5 & \quad 6,7 & $N+1:2N$  \\
        \hline
        4 & \E6{-12} & $-0.946256$ & $-2.1663$ & n/a & $-397.$  \\
        8 & \E{-3}{-12} & $-0.996073$&$ -3.7850$&$ -7.121$&$
        -1588.$ \\
        16 & \E{-1}{-10} & $-0.999750$&$ -3.9843$&$ -8.832$&$
        -6355.$  \\
        32 & 0 & $-0.999986$&$ -3.9990$&$ -8.989$&$ -25421.$  \\
        \hline
    \end{tabular}
\end{table*}

Similarly, consistency errors are~\Ord{H^{2p}} for \(2p-1\) stencil width for general \pde{}s.

Here we only discuss the case when the microscale simulator needs field values~\(u\) on the patch edges.
Analogous formula successfully interpolate derivative values~\(u_x\) to the patch edges if needed by the microscale simulator \cite[]{Roberts04d}, or Robin conditions, or two-point conditions \cite[]{Roberts06d}.

\subsection{and with dynamical systems support}
\label{sswdss}

This section describes one way to provide theoretical support for the patch scheme in its macroscale modelling of nonlinear microscale systems.
My trick is to embed the patch scheme in a one parameter family of schemes. 
The introduced parameter~\(\gamma\) controls the strength of the inter-patch coupling: when \(\gamma=1\) the patches are fully coupled; when \(\gamma=0\) the patches are isolated from each other.
We use a theory based at \(\gamma=0\) to access results that hold for full coupling, \(\gamma=1\).

Introduce the parameter~\(\gamma\) into the inter-patch coupling~\eqref{eqsdoper} so it takes the modified form
\begin{subequations}\label{eqsdopec}%
\begin{align}
    E_\gamma^{\pm r}
    &:=1 +\gamma r(\pm\mu\delta +\rat12r\delta^2 )
    \label{eqdopeca}\\&\quad{}
    +\gamma^2r(r^2-1)(\pm \rat1{3!}\mu\delta^3
    +\rat1{4!} r\delta^4)
    \label{eqdopecb}\\&\quad{}
    +\gamma^3r(r^2-1)(r^2-4)(\pm\rat1{5!}\mu\delta^5
    +\rat1{6!}r\delta^6)
    \label{eqdopecc}\\&\quad{}
    +\gamma^4r(r^2-1)(r^2-4)(r^2-9)(\pm \rat1{7!}\mu\delta^7
    +\rat1{8!}r\delta^8)
    \label{eqdopecd}\\&\quad{}
    +\Ord{\gamma^5}.
\end{align}
\end{subequations}
Then, instead of~\eqref{eqcicpdev}, the patch-edge values are determined by
\begin{equation}
u_j(\pm r,t)=E_\gamma^{\pm r}u_j(0,t), 
\label{eqcicpdec}
\end{equation}
for scale ratio \(r=h/H\)\,.
Observe that terms multiplied by~\(\gamma^1\) flag a patch communicating with its nearest neighbours, those terms multiplied by~\(\gamma^2\) flag communication out to next-nearest neighbours, and so on.
\footnote{That is, parameter~\(\gamma\) empowers us to order the system's interactions between all the patches (maybe analogous with Feynman diagrams in physics).}
Hence, analysing to asymptotic error~\Ord{\gamma^{p+1}} means that a patch communicates with \(p\)~patches to each side of itself.

Incidentally, when \(r=1\) the patches overlap and empowers algebraic novel and accurate discretisations based upon the \pde\ telling us the sub-grid structures (\cref{algmd}).  Then \(E^{\pm1}_\gamma=1+\gamma(E^{\pm1}-1)\), the so-called \emph{holistic discretisation}. 

To get the flavour of the theoretical support this patch coupling engenders for our multiscale computational of patches and gaps, let's consider this \(\gamma\)-parametrised scheme applied to Burgers' \pde\ \(u_t=u_{xx}-uu_x\) with say \(N\)~patches.
\begin{enumerate}
\pgfplotsset{ colormap={allred}{rgb=(0.7,0,0) rgb=(1,0,0)}}
\def\smi#1{\marginpar{
\begin{tikzpicture}
\begin{axis}[footnotesize,font=\footnotesize
    ,xmin=0,ymin=0,zmin=-1,xmax=1.1,ymax=1.2,zmax=1.1, ticks=none
    ,axis lines=middle, view={110}{25}
    ,xlabel={$\gamma$},ylabel={$\vec U,V_H$},zlabel={$\{u_j(x)-U_j\},W_H$}]
\addplot3+[no marks,very thick,forget plot,domain=0:1.1]({0},{x},{0}) 
    node[below left] {equilibria};
\ifnum#1>0
\addplot3+[-stealth,samples=4,no marks,thin,domain=0:1.1,y domain=-1:1
    ,quiver={u=0,v=0,w=-z,scale arrows=0.5,every arrow/.append style={thin}}] 
    ({0},{x},{y});
\fi
\ifnum#1>1
\addplot3+[surf,no marks,domain=0:1,y domain=0:1.1,opacity=0.5
    ,samples=7,colormap name=allred] {y*x^2/2};
\fi
\end{axis}
\end{tikzpicture}
}}%
\item \emph{Equilibria?} When \(\gamma=0\) each patch is isolated, and so \(u_j={}\)constant independently in each patch is an \(N\)D subspace of equilibria (cf. DP's finite element space~\(V_H\)); that is, \(u_j=U_j\) for \(j=1,\ldots,N\)\,---see schematic marginal plot.
\smi0

As it is easiest we analyse about the equilibria \(u_j=0\), but in principal we could generalise to being global in~\(\{U_j\}\).  The analysis is local in coupling parameter~\(\gamma\), but evidence indicates the locale often extends to include \(\gamma=1\)\,---the case of interest.

\item \emph{Linearisation:} small perturbations to the equilibrium then satisfy, upon changing to \(\xi=(x-X_j)/H\), the \pde\ \(\hat u_t=\frac1{H^2}\hat u_{\xi\xi}\) such that, from~\eqref{eqcicpdec} with \(\gamma=0\), \(\hat u_j(\pm r,t)=\hat u_j(0,t)\).
As coupling \(\gamma=0\), each patch is isolated and the analysis here is the same over all patches.
As in \cref{ssndop}, seek solutions via separation of variables, \(\hat u_j=e^{\lambda t}v(\xi)\), and find eigenvalues and eigenfunctions are
\begin{align*}&
\lambda_k=-\pi^2k^2/h^2,
\\&
v_k(\xi)=\begin{cases}
\sin(k\pi\xi/r)&k=1,2,3,\ldots\,,\\
\cos(k\pi\xi/r)&k=0,2,4,\ldots\,,
\end{cases}
\end{align*}
as well as some generalised eigenfunctions for \(k=2,4,6,\ldots\)\,---2nd schematic marginal plot.
\smi1

\item \emph{Emergence theory?} Because there is a zero eigenvalue for each of \(N\)~patches, with all other eigenvalues\({}\leq-\pi^2/h^2\), centre manifold theory \cite[e.g.]{Carr81, Haragus2011} asserts the \pde\ on \(N\)~patches coupled by~\eqref{eqcicpdec} possesses an \((N+1)\)D slow centre manifold~\(u_j=u_j(\xi,\Uv,\gamma)\) on which the system evolves \(\dot U_{j}=g_j(\Uv,\gamma)\), for some domain of finite~\(\gamma\). 
Further, the solutions on this slow centre manifold attract \emph{all} nearby ones roughly as~\(\exp(-\pi^2t/h^2)\).
That is, on a cross-patch diffusion time all solutions of the multi-patch system approach the dynamics of a slow macroscale \(N\)D~system---3rd schematic marginal plot with red centre manifold.
\smi2

This slow centre manifold system turns out to be an appropriate macroscale system at full coupling \(\gamma=1\).

\item \emph{Construction:} computer algebra handles the tedious details.
It eventuates that on the slow centre manifold, the \(j\)th~patch has sub-patch field (cf. DP's \(V_H\oplus W_H\))
\begin{equation*}
u_j=U_j +\gamma
    \left( \xi\mu\delta +\half\xi^2\delta^2 \right)U_j
    +\Ord{\gamma^2+|\Uv|^2}.
\end{equation*}
Here, this sub-patch field is the classic parabola formed from local first and second derivative estimates.
Higher-order linear terms are likewise.
However, microscale heterogeneity, non-linearity or odd derivatives in the \pde\ (\cref{algmd}) generate non-classic sub-patch structures that represent non-trivial out-of-equilibrium structures on the sub-patch microscale \cite[\S3, e.g.]{Roberts00a}.

The evolution on the slow centre manifold is
\begin{equation*}
\dot U_j =
    \frac1{H^2}\left[ \gamma\delta^2 -\rat1{12}\gamma^2\delta^4
    +\rat1{90}\gamma^3\delta^6 \right] U_j
    +\Ord{\gamma^4+|\Uv|^2}.
\end{equation*}
Choosing to truncate with coupling errors~\Ord{\gamma^{p+1}} and then evaluating at full coupling \(\gamma=1\) gives classic spatial discretisations of the diffusion \pde, with errors classically recorded as~\Ord{H^{2p}}.
That is, \emph{the slow centre manifold emergent dynamics of the coupled small patches is precisely a sound model of the macroscale dynamics.}

The same conclusion holds when the nonlinear terms are explored, and also holds in two space dimensions \cite[]{Roberts2011a}.

\end{enumerate}

\emph{We conjecture the patch-scheme is similarly good for scenarios where we do not know the macroscale closure.}

\subsection{Open problems}
\begin{itemize}
\item For periodic patches: investigate other controllers?  their optimal control? extend supporting analysis to a useful range of macroscale \pde{}s? including stochastic? develop coupling in multiple directions, not just 1D?  develop boundary patches?

\item A trendy activity is to get Deep Neural Networks (!) to `learn' a macroscale closure from microscale simulations (\cref{algmd}): compare such \textsc{dnn}s with algebraic closures \cite[e.g., cf.][]{BarSinai2018, Roberts00a}.

\item \cite{Roberts08l, Roberts2011a} provide theoretical support for the patch scheme on a regular grid in 2D: it should be straightforward to extend to more space dimensions; and may be challenging to extend to an unstructured grid of patches.

\item We are currently working on developing moving patches with the aim of capturing  shocks by resolving them on a microscale patch and without assuming any Rankine--Hugoniot conditions.

\item Develop the toolbox to patch functions for non-periodic macroscale boundary conditions, to higher-D, to effectively parallelise, \text{and so on}.
\end{itemize}

\section{Macroscale computation of microscale spatial complexity}
\label{smcmsc}
\localtableofcontents

Recall that \cref{smcms} discussed the amazing automatic homogenisation of one patch period-two heterogeneous diffusivity.
Here we discuss multiple patches and the corresponding automatic macroscale homogenisation of the patch scheme.

We are developing a suite of \script\ functions to empower users to take advantage of the patch scheme and other multiscale techniques.
Download the current version from GitHub. 
\footnote{\protect\url{https://github.com/uoa1184615/EquationFreeGit.git}}

Many of the main functions, if invoked with \emph{no} arguments, will execute a basic example.  
For example, executing \verb|configPatches1()| draws \cref{fig:burg3} arising by simulating Burgers' \pde\ within eight patches.  

The user manual, \verb|eqnFreeUserMan-newest.pdf|, is in the main folder and should suffice for most users.
\footnote{\texttt{Doc/eqnFreeDevMan.pdf} documents full details of the toolbox functions.}

\subsection{Couple patches of microscale heterogeneous diffusion for macroscale accuracy}
\label{sscpmhdma}

The script \verb|homogenisationExample| simulates the basic homogenisation introduced by \cref{sammt}, but now on multiple patches, albeit still in~1D.
Recall the remarkable result that if configured so that the patch half-width is an integral multiple of the microscale period, then the patch scheme simulates the exact macroscale homogenisation \cite[\S5.2]{Bunder2013b}.

The overall plan of the code is similar to that discussed in \cref{sspfeft}, but microscale details are different.
The user has to drive two functions in the toolbox: \begin{itemize}
\item \verb|configPatches1()| configures the arrangement of patches and sub-patch microscale lattice in the domain;   
\item \verb|patchSmooth1()| computes the patch edge-values so that a user's function computes a time-step\slash derivative of the sub-patch structure.
\end{itemize}
The overall plan is the following:
\begin{enumerate} 
\item invoke \verb|configPatches1()| and other initialisation
\item user time loop/integration, e.g.~\verb|ode15s|
\begin{enumerate}
\item invoke \verb|patchSmooth1()|
\begin{enumerate}
\item \verb|patchEdgeInt1()| computes patch edge values
\item user function \verb|heteroDiff()| gives time derivatives
\end{enumerate}
\end{enumerate}
\item process results
\end{enumerate}
The code initialises the heterogeneous diffusion coefficients and periodicity
\begin{matlab}
mPeriod = 3
cDiff = exp(randn(mPeriod,1))
\end{matlab}
Then it invokes \verb|configPatches1()| to configure a \(2\pi\)-periodic domain with nine patches coupled with fourth-order interpolation.  Each patch is chosen to be of size ratio~\(0.2\) and to contain precisely two periods of the microscale heterogeneity.
\begin{matlab}
global patches
nPatch = 9
ratio = 0.2
nSubP = 2*mPeriod+1
Len = 2*pi;
ordCC = 4;
configPatches1(@heteroDiff,[0 Len],nan ...
    ,nPatch,ordCC,ratio,nSubP);
\end{matlab}
The code takes advantage of the patch struct~\verb|patches| to communicate the specific heterogeneous coefficients, identical for each patch (could they be different?),  
\begin{matlab}
patches.c = repmat(cDiff ...
    ,(nSubP-1)/mPeriod,1);
\end{matlab}
to the user's microscale function \verb|heteroDiff()|.
\begin{matlab}
function ut = heteroDiff(t,u,x)
  global patches
  dx = diff(x(2:3)); % space step
  i = 2:size(u,1)-1; % interior patch
  ut = nan(size(u)); % preallocate
  ut(i,:)=diff(patches.c.*diff(u))/dx^2; 
end
\end{matlab}

The example script then integrates in time from some initial condition using \verb|ode15s()|
\begin{matlab}
u0 = sin(patches.x) ...
    +0.4*randn(nSubP,nPatch);
[ts,ucts] = ode15s(@patchSmooth1 ...
    , [0 2/cHomo], u0(:));
\end{matlab}
To visualise the results, plot as before
\begin{matlab}
xs = patches.x;  xs([1 end],:) = nan;
mesh(ts,xs(:),ucts'),  view(60,40)
xlabel('time t'), ylabel('space x')
zlabel('u(x,t)')
\end{matlab}

\paragraph{Accuracy?}
\cite{Bunder2013b} [\S5.2] analyse the scenario and assure us that the macroscale predictions are correct.
We may verify here by computing the Jacobian of the patch scheme, and then the small magnitude eigenvalues correspond to the macroscale modes.
I wrote the script \verb|homogenisationAccuracy| to do this.
It normalises the heterogeneous diffusions so their harmonic average is one and so the homogenised \pde\ is \(U_t=1\cdot U_{xx}\).
Then on the \(2\pi\)~domain the macroscale eigenvalues should be~\(-n^2\) for integer~\(n\).
The script uses spectral interpolation to `eliminate' interpolation errors.
Executing the script gives answers such as the following.
\needspace{12\baselineskip}
\begin{multicols}2
\begin{verbatim}
mPeriod =
     3
cDiff =
        5.887
        21.52
      0.35924
nPatch =
     9
ratio =
          0.1
nSubP =
     7
lamFast =
      -5440.7
lam0 =
      -7.7342e-12
      -0.9997
      -0.9997
      -3.9947
      -3.9947
      -8.9731
      -8.9731
      -15.915
      -15.915
\end{verbatim}
\end{multicols}
These macroscale eigenvalues are close to~\(-n^2\). 
The differences are not errors, the differences are due to the small but finite size of the microscale discrete lattice causing higher-order terms in the effective macroscale \pde:  
it should be \(U_t=1\cdot U_{xx}+?d^2U_{xxxx}+\cdots\).
Reduce the ratio to lessen higher-order effects from the microscale.

\begin{figure*}
\centering
\caption{\label{fig:HomogenisationMicro}cross-eyed stereo
pair of the field~\(u(x,t)\) during each of the microscale
bursts used in the projective integration of heterogeneous
diffusion.}
\includegraphics[scale=0.99]{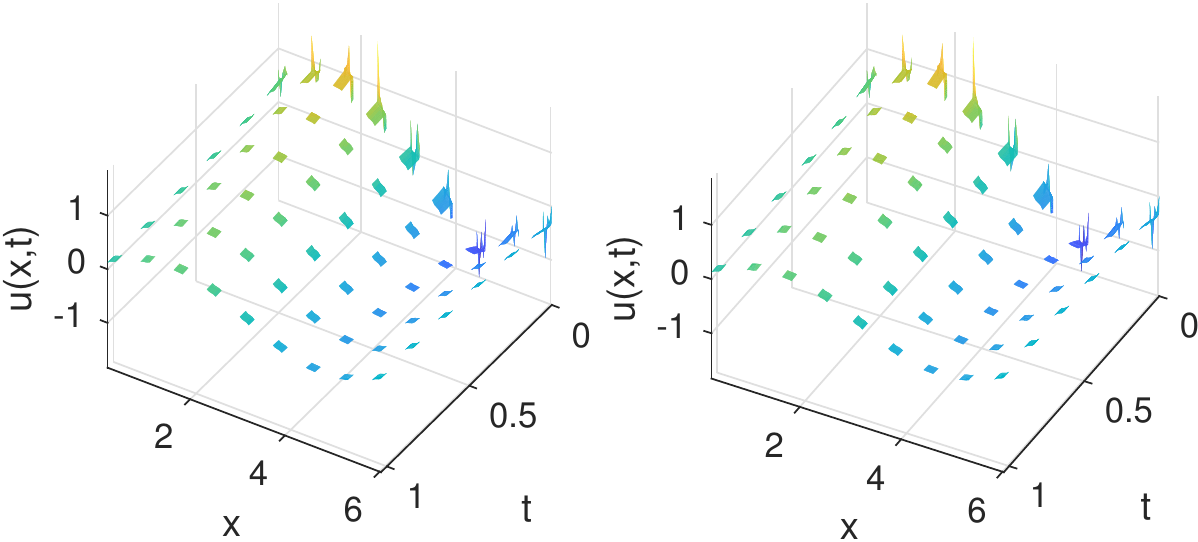}
\end{figure*}
\paragraph{Projective integration}
There is a large spectral gap between these macroscale modes and the leading eigenvalue of the sub-patch modes, by a factor of roughly a thousand.
These multiscale schemes are typically extremely stiff so prefer projective integration (\cref{spicosbt}). 
The script \verb|homogenisationExample| proceeds to additionally use projective integration in time, and \cref{fig:HomogenisationMicro} illustrates that the resulting simulation only computes on patches in space-time.

\subsection{Avoid buffers in general homogenisation}

So far we have addressed the scenario where we know the periodicity of the microscale heterogeneity.  What if we do not? or if the microscale is random?  These were explored by \cite{Bunder2013b} [\S5.3].

\begin{figure*}
\centering
\caption{\label{fig:onepatch}The \emph{core} region in the centre of this patch, $|i|\leq c$\,, is used for the macroscale amplitude. 
The two outlined \emph{action} regions on the ends of the patch are used for the coupling. 
The so-called buffers in between have width $b=n-c$ \protect\cite[Fig.~4]{Bunder2013b}.}
\includegraphics{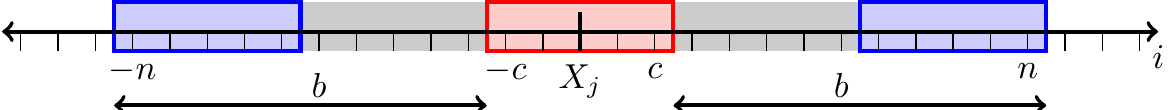}
\end{figure*}
\cref{fig:onepatch} illustrates a more general patch of size \(2n+1\) lattice points.
Instead of defining the patch macroscale value as the centre patch value, we define it to be the average over a \emph{core region} of width \(2c+1\) at the centre of a patch.
Since the macroscale is an average, the coupling conditions also need to be phrased in terms of corresponding averages.
Hence patches are coupled by defining \emph{action regions} of width \(2c+1\) at the edges of each patch, and then requiring that the average over an action region be the interpolated macroscale averages.
This forms a more general patch scheme.

\begin{SCfigure*}
\centering
\caption{\label{fig:structure}Microscale structure of a single patch with microscale period $K=3$ and different core half-widths~$c$: the core and action regions are shaded and outlined \protect\cite[Fig.~7]{Bunder2013b}.}
\includegraphics[width=0.7\textwidth, 
height=\textheight,  keepaspectratio]{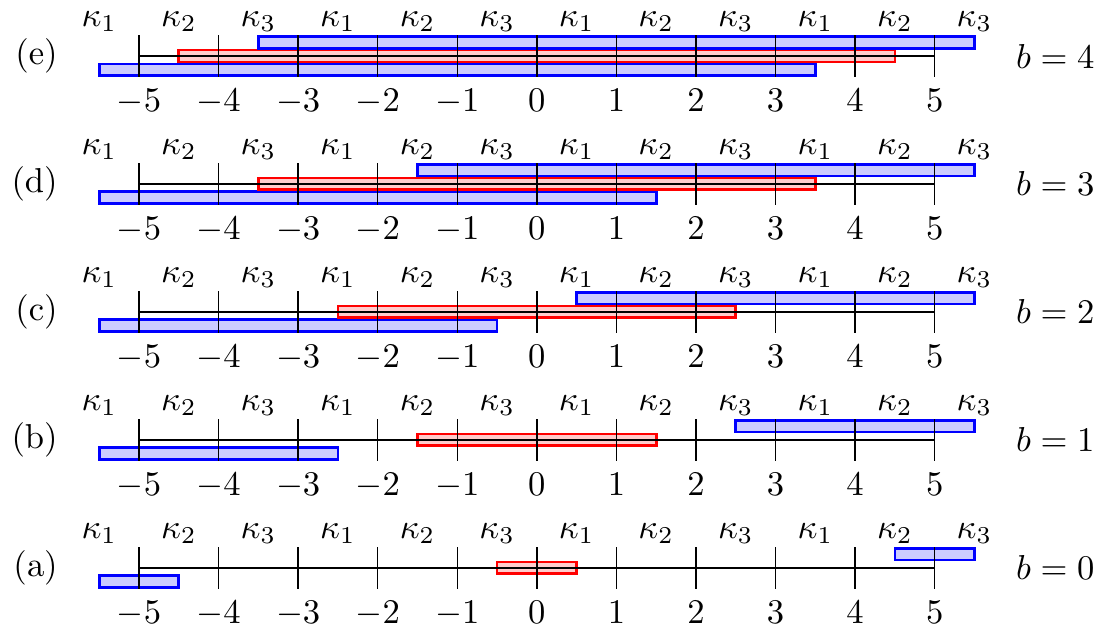}
\end{SCfigure*}
We explored this more general patch scheme in simulating the macroscale dynamics of microscale heterogeneous diffusion for cases where we do not `know' the microscale periodicity.
\cref{fig:structure} shows we varied the size of the core\slash action regions, assessed the error in the macroscale predictions over a range of microscale periodicities.
Be particularly interested in the commonly espoused plausible idea that there should be a sizeable \emph{buffer} (grey) between the action regions and the core region to allow the sub-patch \text{solution to `heal'.}

\begin{SCfigure*}
\centering
\caption{\label{fig:allk0}Coefficient relative errors averaged over microscale periods $2\leq K\leq 12$\,, including \(K>n\)\,, versus the relative core half-width $0\leq c/(n-1)\leq 1$ for patch sizes~$n$. 
The error~$|\rho_0|$ is minimised when $c\approx 0.4n$ \protect\cite[Fig.~9]{Bunder2013b}.}
\includegraphics[width=0.7\textwidth, 
height=\textheight,  keepaspectratio]{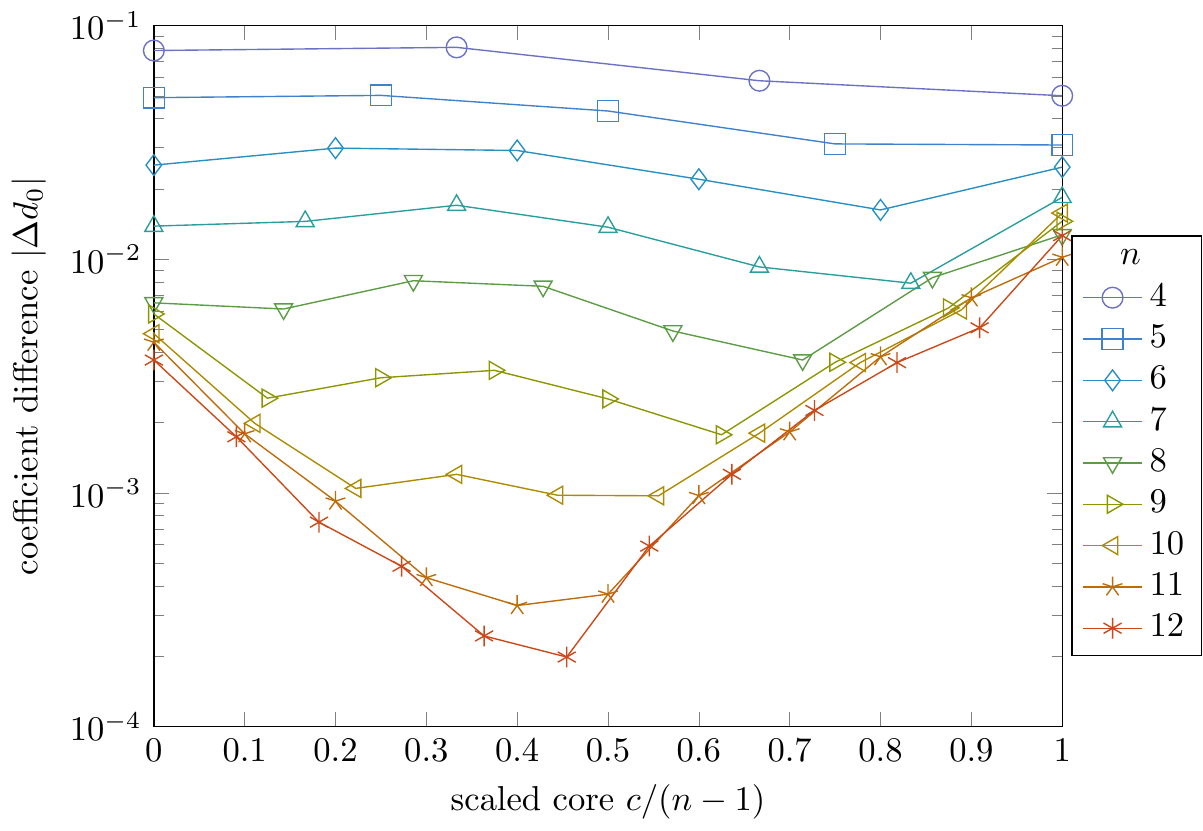}
\end{SCfigure*}
\cref{fig:allk0} shows the log-error as a function of the core size relative to the patch size.
A remarkable feature of \cref{fig:allk0} that the error is minimised for core \(c\approx 0.4 n\).
That is, for minimum error in macroscale predictions make the core and action regions overlap slightly!  
Making the core and action regions abut without overlapping is also reasonable if you prefer as \(\frac13\approx 0.4\).
It appears that we never need buffers.

\subsection{How do communication delays affect such simulations?}

Recall we envisage that in large problems the spatial patches will be distributed across many processors in parallel.  
In that scenario much of the microscale computation may be done efficiently on each processor.
However, then the inter-patch coupling incurs real-time expensive inter-processor communication.
In the scheme described so far, such inter-patch\slash processor communication occurs \emph{every} microscale time step. 
Such communication would slow the computation severely.

\begin{SCfigure*}
\centering
\caption{\label{figmesofig}Three patches are the shaded regions centred on macroscale lattice points~$X_{i-1}$\,, $X_i$ and $X_{i+1}$\,.  
At times \(\delta t\)~apart the coupling information from neighbouring patches is updated \protect\cite[Fig.~3]{Bunder2015a}.}
\includegraphics[width=0.75\textwidth, 
height=\textheight,  keepaspectratio]{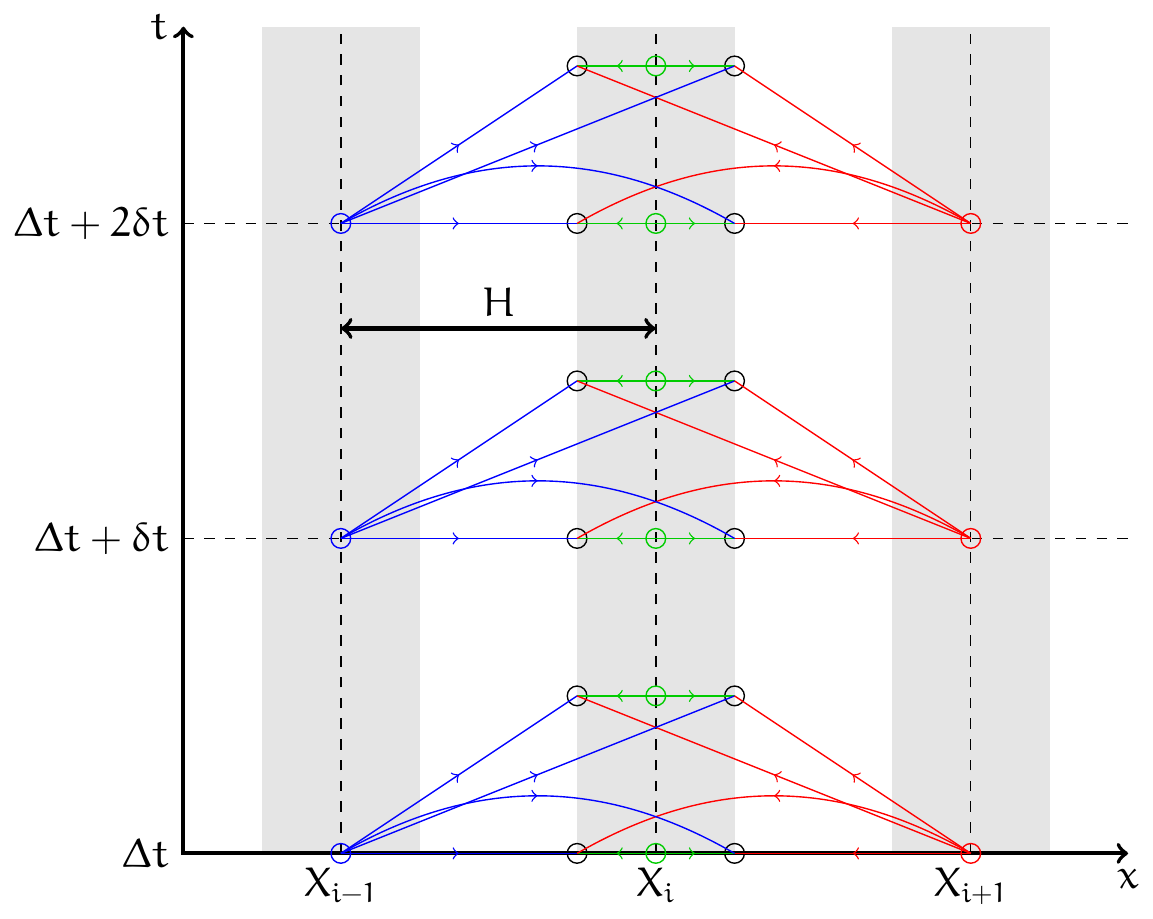}
\end{SCfigure*}%
Alternatively, maybe we could communicate the inter-patch coupling less often?
\cite{Bunder2015a} explored the effect on accuracy of communicating coupling on an intermediate time meso-scale---longer than micro-times, and shorter than macro-times.
\cref{figmesofig} illustrates \text{the idea.}

For simple diffusion on a lattice \(\dot u_j(t)=u_{j+1}(t)-2u_j(t)+u_{j-1}(t)\), we could obtain complicated analytic formulas across many scenarios.
To help fill-in the communication `gap' from time~\(t\) to time~\(t+\delta t\) we considered communicating not only the neighbouring macroscale values, but also their first~\((Q-1)\) time derivatives. 
During a time~\(\delta t\) coupling errors penetrate into the core of each patch from the edge regions.  The aim is to keep small the errors in the patch-core.

\begin{SCfigure*}
\centering
\caption{\label{figpenetrateA20}the upper bound of components of the remainder\slash error~$R_{j\max}$ in a patch with patch half-width $n=20$\,, for mesoscale time~$\delta t=0.5$  and $Q=1,3,5,7$ \protect\cite[Fig.~10]{Bunder2015a}.}
\includegraphics[width=0.75\textwidth, 
height=\textheight,  keepaspectratio]{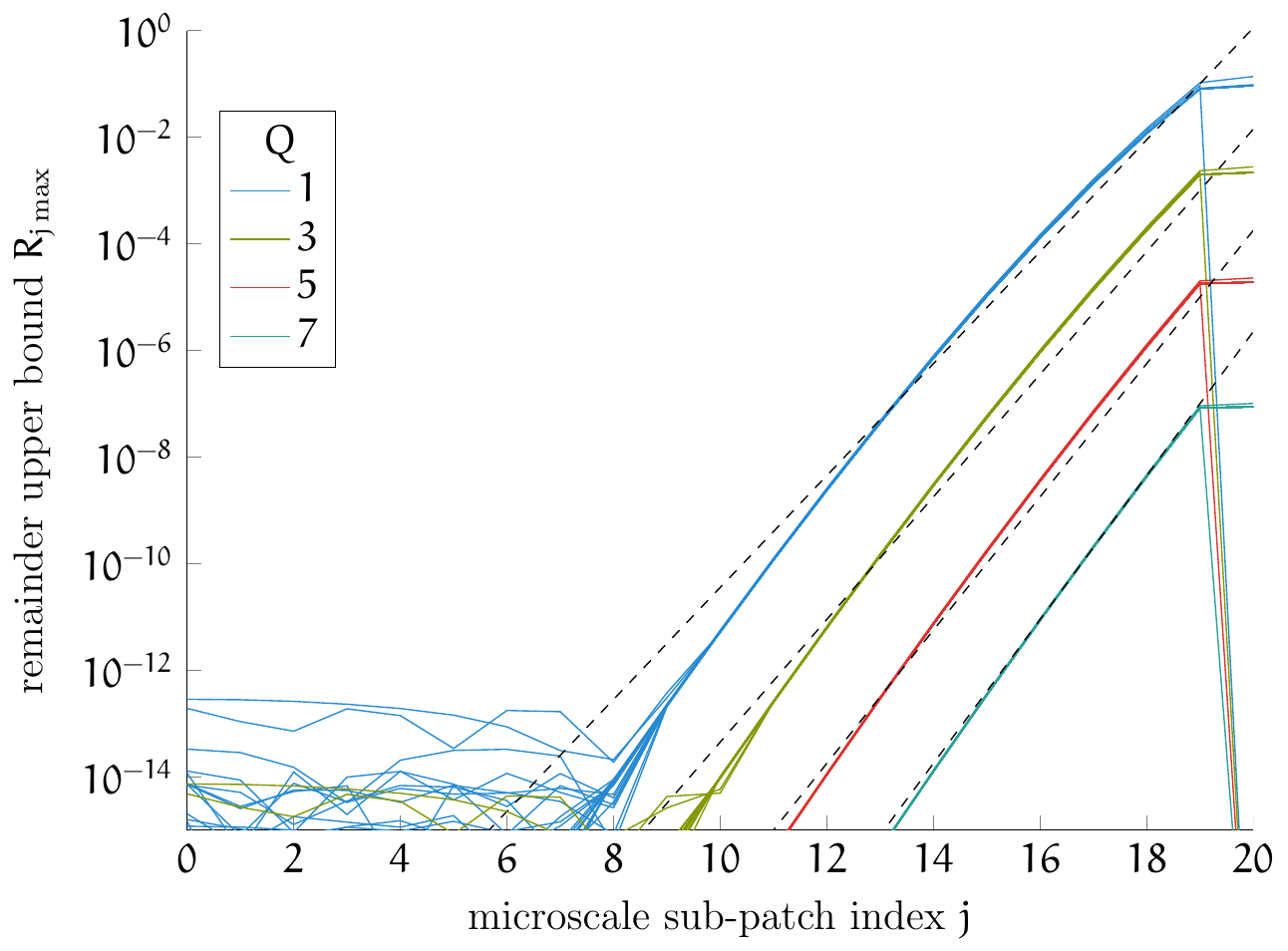}
\end{SCfigure*}%
\cref{figpenetrateA20} shows that by using larger patches, here half-size \(20\)~lattice points, we protect the patch-core from the lack of updates to the patch edge-values.
Here the error in the core is at the round-off level.
Increasing the number~\(Q\) of time derivatives communicated decreases the error.
\cref{figpenetrateA20} shows the results for all cases where a core-average over \(2a+1\) points in the patch is the macroscale quantity communicated.

\begin{SCfigure*}
\centering
\caption{\label{figerror1}the error in the core,~$E_{\max}$, for $Q=1$\,, only function values communicated, a range of  mesoscale time steps~$\delta t$ and over several patch half-widths~$n-a$ \protect\cite[Fig.~11]{Bunder2015a}.}
\includegraphics[width=0.75\textwidth, 
height=\textheight,  keepaspectratio]{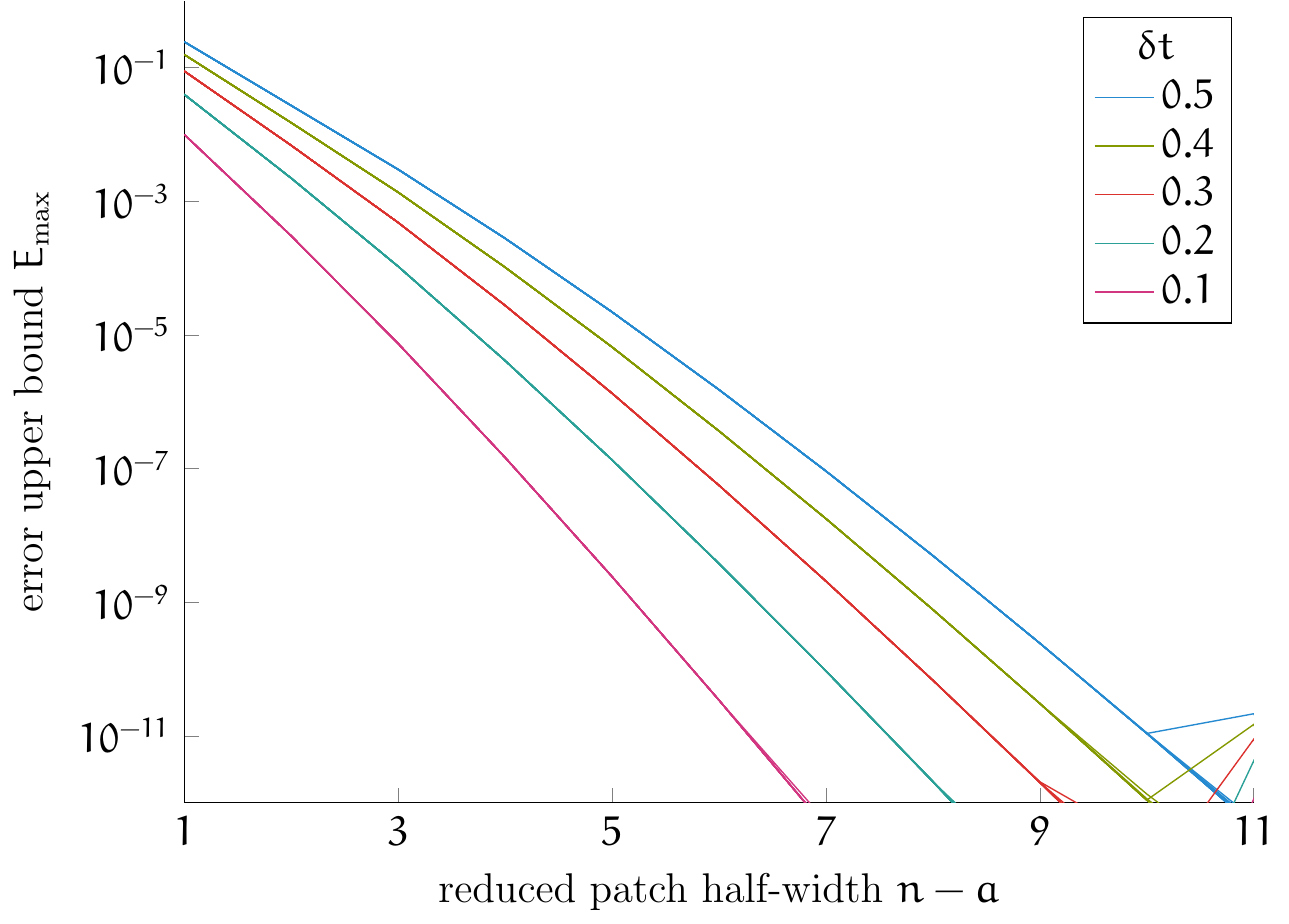}
\end{SCfigure*}
\cref{figerror1} shows the effect of varying the meso-scale time of communication,~\(\delta t\).  As expected the errors increase with~\(\delta t\), and decrease with increasing patch width.
Similar behaviour was seen in 2D spatial \text{patch simulations.}

Limiting inter-patch\slash processor communication to meso-times may speed up large-scale simulations by factors of \(10\)--\(1000\) \cite[\S6.2]{Bunder2015a}.
But in limited communication we need \text{some `buffers'.}

%\begin{activity}
%What is a toy one-patch problem that illustrates key results here?
%\end{activity}

\subsection{Algebra `learns' good macroscale discretisations}
\label{algmd}

For efficient computation within, we want small patches.
But when algebra resolves the sub-patch structure then the patches can be any size.
Patches that overlap, \(r=1\), appear most appealing---the so-called \emph{holistic discretisation} \cite[]{Roberts00a}.

This approach provides a framework for systematically constructing sub-grid structures determined by the problem---like DP, but non-`variational'---in contrast to \textsc{fd/fe/fv} imposition.

We explore the macroscale dynamics predicted by overlapping patches.
We find formulas for sub-patch structures and macroscale evolution that are cognate to trendy efforts that analogously try to `discover' macroscale closures by invoking Deep Neural Networks on extensive numerical simulations \cite[e.g.]{BarSinai2018}.

Here let's explore~\(u(x,t)\) governed by the linear advection-diffusion \pde
\begin{equation}
u_t=-c\D xu+\DD xu\,.
\label{eqadpde}
\end{equation}
Define a macroscale grid~\(\{X_j\}\) of equi-spacing~\(H\).
Let the \(j\)th~patch be \(|x-X_j|\leq H\) so it stretches from~\(X_{j-1}\) to~\(X_{j+1}\), and the field on the \(j\)th~patch be~\(u_j(\xi,t)\) for sub-patch variable \(\xi=(x-X_j)/H\).

\emph{2D space?}  There are several ways to generalise to 2D \cite[]{Roberts08l, Roberts2011a}.  The overlap means a patch shares half of itself with each of four nearest neighbours.  Hence the overlap empowers `channels' to be continued or otherwise correctly from one patch to the next---in principle (as required by BL in so-called `skins').

In 1D the coupling between patches is particularly simple as the edge of each patch is the neighbouring centre-patch.
Then the embedded classic interpolation~\eqref{eqcicpdec} becomes simply
\begin{equation}
u_j(\pm 1,t)=(1-\gamma)u_j(0,t)+\gamma u_{j\pm1}(0,t)
\label{eqhdpdec}
\end{equation}
in terms of coupling parameter~\(\gamma\), and local space variable~\(\xi\).

As in patches, the parameter~\(\gamma\) controls the \emph{locality} in space of a hierarchy of approximations: working to `errors'~\(\Ord{\gamma^p+1}\) means a stencil width \text{of \(2p+1\).}

Our quest is to find that the macroscale dynamics of the patch system~\eqref{eqadpde,eqhdpdec} is good.
\begin{enumerate}
\item Equilibria?  When the coupling \(\gamma=0\), the patches are isolated, and equilibria are that \(u_j(\xi,t)=\)constant independently in each patch.
Defining the amplitude to be the centre-patch value, \(U_j(t)=u_j(0,t)\), these constants are~\(U_j\).
\marginpar{The schematic marginal plots of \cref{sswdss} apply here also.}

\item Linearise: each patch is the same with linearised problem \(H^2u_{j,t}=-cHu_{j,\xi}+u_{j,\xi\xi}\) such that \(u_j(\pm 1,t)=u_j(0,t)\).
Each patch has a zero eigenvalue, and an infinite number of negative eigenvalues: maybe \(\lambda\leq -\pi^2/H^2-c^2/4\). 

\item Hence, solutions exponentially quickly approach a slow centre manifold\slash subspace that is parametrised by~\(\{U_j\}\) and~\(\gamma\) \cite[e.g.,][]{Carr81, Haragus2011}.

\item Construct?  
\footnote{Although the problem is linear, the \emph{construction of a model is nonlinear} through the chain rule that \(u_{t}=(\D Uu)U_t=(\D Uu)g\) and we must determine both~\(u(U)\) and~\(g(U)\).}
Seek sub-patch fields \(u_j=U_j+\gamma v_j(\xi,\Uv)+\Ord{\gamma^2}\) such that \(\dot U_j=\gamma g_j(\Uv)+\Ord{\gamma^2}\).
Substitute into coupled \pde~\eqref{eqadpde,eqhdpdec} and the \(\gamma\)-terms require
\begin{align*}&
H^2g_j=-cHv_{j,\xi}+v_{j,\xi\xi}
\\&\text{s.t. } 
v_j(\pm1,t)-v_j(0,t)=U_{j\pm1}-U_j\,.
\end{align*}
Although tedious, the solution is straightforward \cite[\S3]{Roberts00a}:
\begin{subequations}\label{eqs:}%
\def\heh{(cH/2)}\def\ex{{cH\xi}}\def\eh{{cH}}%
\begin{align}
    v_j & =  \left[ \frac{e^\ex-1}{4\sinh^2\heh} 
    -\frac{\cosh\heh}{2\sinh\heh}\xi \right]\delta^2U_j 
    \nonumber\\&\quad{}
    +\xi{\mu\delta} U_j\,,
    \label{eq:v1}  \\
    g_j & =  -\frac{\mu\delta}{H} U_j
    +\nu_1(cH)\frac{\delta^2}{H^2} U_j \,,
    \label{eq:g1}  \\
    \nu_1 & =  \frac{\eh\cosh\heh}{2\sinh\heh}\,.
    \label{eq:nu1}
\end{align}
\end{subequations}
Higher orders in coupling~\(\gamma\) are best left for computer algebra.

\item Interpret?  Set coupling parameter \(\gamma=1\) to predict at full coupling.
\begin{itemize}
\item Expression~\eqref{eq:v1}, added to~\(U_j\), gives the sub-patch structures `learnt' by the algebra from the \pde.  
These structures are \emph{not imposed by us} on the physical problem  (cf.~classic finite element\slash volume\slash differences).
These out-of-equilibrium structures `know' relevant sub-patch dynamics from the \pde.
\item Higher order analysis in coupling~\(\gamma\), say to errors~\Ord{\gamma^{p+1}}, would learn more about the sub-patch structures by accounting for the influence of patches up to \(p\)-distant \cite[\S3]{Roberts00a}.

\begin{figure}
\centering
\caption{\label{fignu1}the enhanced macroscale dissipation~\eqref{eq:nu1}.}
\begin{tikzpicture}
\begin{axis}[small,domain=0.02:5.5,ymin=0,axis lines=middle
,xlabel={$cH$},ylabel={$\nu_1$}]
\addplot+[no marks]{x/2*cosh(x/2)/sinh(x/2)};
\end{axis}
\end{tikzpicture}
\end{figure}
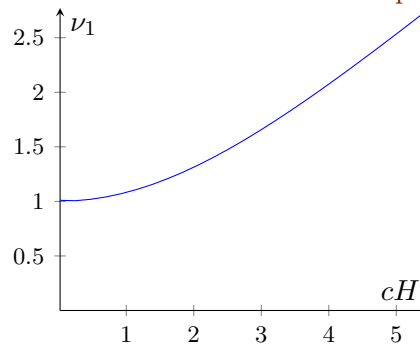
\item The macroscale dynamics are the \ode{}s \(\dot U_j=g_j\) from~\eqref{eq:g1} and give us as good an approximation as one could expect from any three-point stencil.
Importantly, the effective dissipation~\(\nu_1\) (\cref{fignu1}) is enhanced so that as advection~\(c\) increases then \(\dot U_j=g_j(\Uv)\) morphs seamlessly into an upwind discretisation. 
The \ode{}s \(\dot U_j=g_j\) are stable for all~\(cH\)---there is no \textsc{cfl} constraint.
\item The algebraic closure `learnt' here is valid globally in~\Uv, \(c\) and~\(H\): the closure is not limited by the finite extent of the simulations typically underlying any machine learning.
\end{itemize}

\end{enumerate}

\needspace{4\baselineskip}
\emph{Modelling in general?}  \emph{Dynamic} macroscale models are \emph{nonlinear} transformations of the microscale system, here \(\nu_1(cH)\) (\cref{fignu1}), and `harmonic mean' in homogenisation.  So endemic linear arguments necessarily have deficiencies in dynamics.  Consequently, ``the whole is more than the sum of the parts'' because the whole is a nonlinear transform of the parts.

\subsection{Open problems}

\begin{itemize}
\item Ongoing research is exploring patch configurations that automatically homogenise 2D heterogeneous diffusion to high accuracy.
\item Further explore `symmetry' preserving patch\slash holistic coupling conditions \cite[]{Roberts08l} especially for general heterogeneity as discussed by DP \cite[cf.][]{Bunder2013c}, or piecewise linear based coupling that also connects to splines \cite[]{Jarrad2016a}.
\item Convincing people that negative probabilities and negative concentrations are \textsc{ok} \cite[e.g.,][Exercise~5.2]{Roberts2014a}.
\item Develop the toolbox to patch functions coping with stochastic sub-patch structures.
\item Explore meso-time communication in practice, and for advection for which I expect a worse performance.
\item Incorporate tolerance to hardware failure as in massive parallelism, needed for the largest problems, it is likely that one of the millions of \textsc{cpu}s will fail during a computation.
\end{itemize}

\section{Projective integration computes only on small bursts of time}
\label{spicosbt}
\localtableofcontents

\begin{SCfigure*}
\centering
\caption{\label{figegPIMM2}projective integration of Michaelis--Menten enzyme kinetics for scale separation parameter \(\epsilon=0.05\).}
\includegraphics[width=0.7\textwidth, 
height=\textheight,  keepaspectratio]{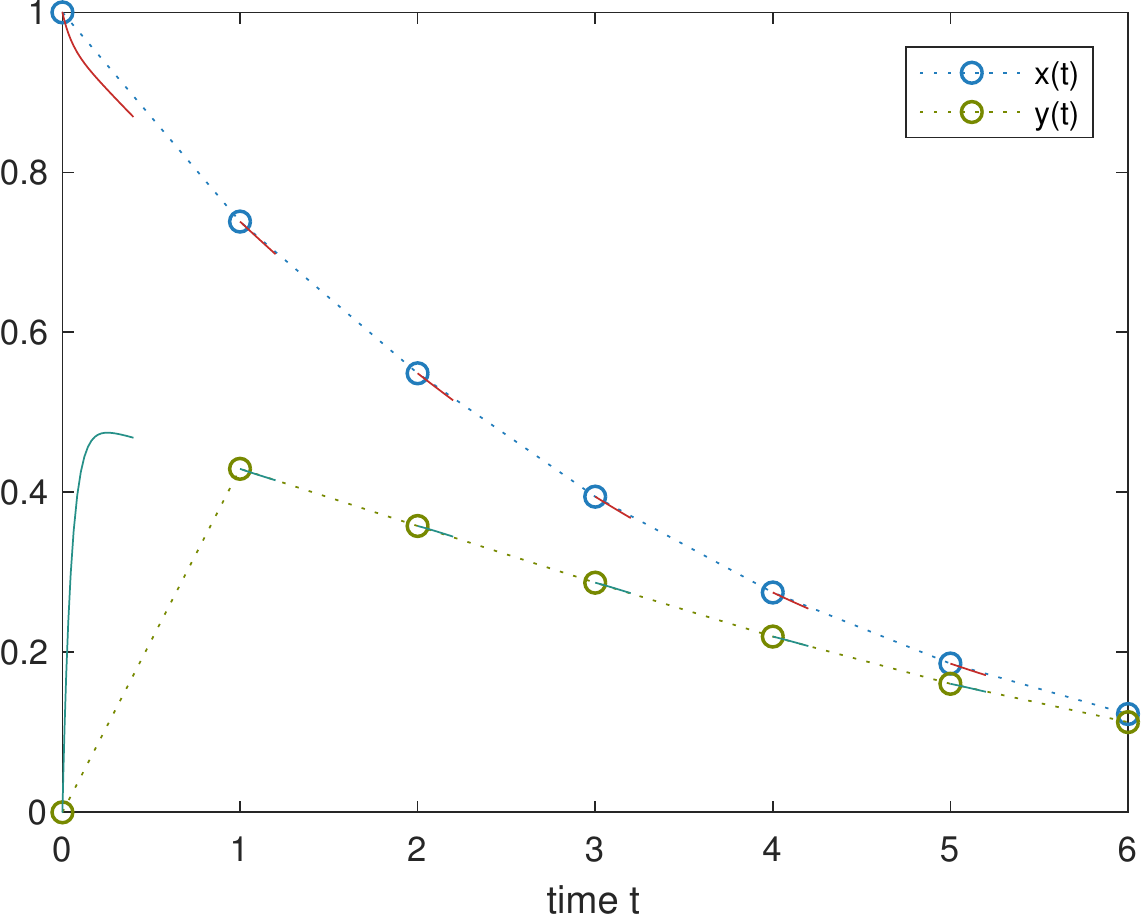}
\end{SCfigure*}%
Historically, projective integration was the first `equation-free' technique \cite[e.g.]{Gear02, Gear02b}.
As shown in the example of \cref{figegPIMM2}, the aim is to use small bursts of microscale simulation (solid lines) to then extrapolate forward in time, over unsimulated time (dotted), to predict the next macroscale value (circles), and then repeat.

Importantly, in this scheme we do \emph{not impose our subjective opinion} of what should be `frozen' macroscale variables.
Some other multiscale methods do subjectively `freeze'.  
Instead projective integration always deals with the full dynamics of the system's complete out-of-equilibrium interactive exchange of information to and fro between micro- and macro-scales.

\subsection{Accuracy and stability for such schemes,}
\label{ssasss}

\begin{figure*}
\centering
\caption{\label{figrowth}for scale ratio \(r=\delta/\Delta=\tfrac19,\tfrac29,\tfrac39\) (blue, red, brown), plot the growth rates~\(G\) as a function of~\(\lambda\Delta\), and compare to (black)~\(e^{\lambda\Delta}\).}
\begin{tikzpicture}
\begin{axis}[small,domain=-2:1,axis lines=middle
,xlabel={$\lambda\Delta$},ylabel={$G$}]
\addplot[no marks,forget plot]{exp(x)};
\foreach \r in {1/9,2/9,3/9} {
  \addplot+[no marks]{exp(\r*x)*(1+x*(1-(\r)))};
  }
\end{axis}
\end{tikzpicture}
\begin{tikzpicture}
\begin{axis}[small,domain=-10:1,axis lines=middle
,xlabel={$\lambda\Delta$},ylabel={$G$}]
\addplot[no marks,forget plot]{exp(x)};
\foreach \r in {1/9,2/9,3/9} {
  \addplot+[no marks]{exp(\r*x)*(1+x*(1-(\r)))};
  }
\end{axis}
\end{tikzpicture}
\end{figure*}
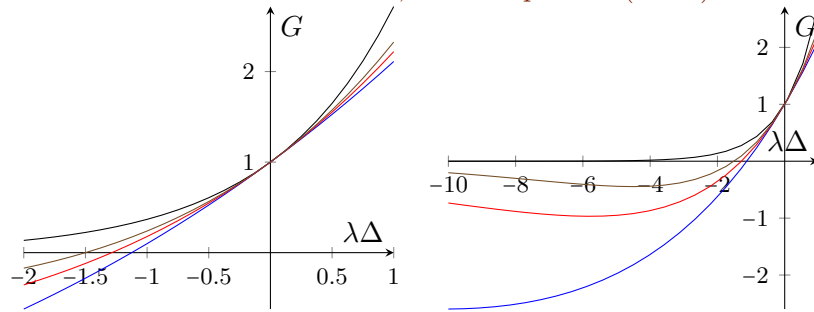

Let's see what happens in the simplest projective integration, akin to Euler method \cite[]{Gear02b}.
You execute a burst of some code and it takes very many microscale steps in time from \(t=0\), say, to arrive at some end-burst time~\(\delta\).
For simplicity, and unknown to you, let's suppose the system is simulating \(u(t)=u_0e^{\lambda t}\).
At the end of the burst,~\(t=\delta\), you get the two results that \(u(\delta)=u_0e^{\lambda\delta}\) and it is changing in time as \(\dot u(\delta)=\lambda u_0e^{\lambda\delta}\).
Extrapolating over time-gap~\((\Delta-\delta)\) then predicts, without any further expensive microscale simulation, that 
\begin{equation*}
u(\Delta)=u(\delta)+(\Delta-\delta)\dot u(\delta)=u_0e^{\lambda\delta}[1+\lambda(\Delta-\delta)].
\end{equation*}
Repeat \(n\)~times to predict \(u(n\Delta)=u_0G^n\) for growth factor,
where \(r=\delta/\Delta\),
\begin{equation}
G=e^{\lambda\delta}[1+\lambda(\Delta-\delta)]=e^{\lambda\Delta r}[1+\lambda\Delta(1-r)].
\label{eq:pigrowth}
\end{equation}
\begin{itemize}
\item \cref{figrowth}(left) compares the growth factor with the exact~\(e^{\lambda\Delta}\) and shows that for accuracy we need small~\(\lambda\Delta\); that is, the macro-time-step must be smallish compared to the time-scale of macroscale interest.
\item \cref{figrowth}(right) shows that for stability the burst needs to be long enough so that \(\delta/\Delta\geq2/9\) (more precisely\({}>0.2178\))---although the bound is much better for a multi-scale system.
\end{itemize}

But in general we must consider large-scale systems, not just one variable.

\paragraph{Systems are locally linear}
Consider a general system of \ode{}s \(\dot\uv=\fv(\uv)\).
Such nonlinear systems are locally linear by Taylor's theorem. 
Suppose \(\uv_*(t)\) is a solution, then nearby solutions   \(\uv(t)=\uv_*(t)+\epsilon {\hat u}(t)\) satisfy
\begin{equation}
\dot {\hat u}=J{\hat u}+\Ord{\epsilon} 
\quad\text{where Jacobian }J:=\left.\D\uv\fv\right|_{\uv_*}.
\label{eqloclin}
\end{equation}
We see this by substituting into the system \(\dot\uv=\fv(\uv)\): the left-hand side is \(\dot\uv_*+\epsilon\dot {\hat u}=\fv(\uv_*)+\epsilon\dot{\hat u}\) whereas by Taylor about~\(\uv_*\) the right-hand side 
\begin{align*}
\fv(\uv_*+\epsilon {\hat u})
&=\fv(\uv_*)+J\epsilon{\hat u}+\Ord{\epsilon^2}.
\end{align*}
Cancelling and dividing gives the linear system~\eqref{eqloclin}.

Recall that, generically, the linear system \(\dot{\hat u}= J{\hat u}\) is fully understood by diagonalisation.  
Generally there exists a linear coordinate change \(\hat u=V\tilde u\), the columns of~\(V\) are eigenvectors, such that \(\dot{\tilde u}=\Lambda\tilde u\) for diagonal matrix of eigenvalues~\(\Lambda\).
By superposition we thus primarily need to consider \(\dot u_i=\lambda_iu_i\) \emph{over all the spectrum} of eigenvalues~\(\lambda_i\) of the Jacobian~\(J=\D\uv\fv\)\,.%
\footnote{Neglecting issues associated with a time varying Jacobian;  for example, see the Marcus--Yamabe system.}

\paragraph{Fast-slow multiscale systems}
Now in a multiscale system, including those with microscale heterogeneity, the spectrum divides into two `clusters' (e.g., \cref{ssndop,sswdss,sscpmhdma}): 
the small eigenvalues of the interesting macroscale, say \(|\Re\lambda_i|\leq\alpha\); and the large negative eigenvalues of the microscale quasi-equilibration, say \(\Re\lambda_i\leq-\beta<\text{gap}<-\alpha\).
In this scenario, our earlier simple analysis~\eqref{eq:pigrowth} indicates the following:
\begin{itemize}
\item choose macrostep~\(\Delta\) such that \(\alpha\Delta\) is small enough for desired accuracy;
\item to ensure the growth~\(G\) satisfies \(|G|<1\) for every \emph{microscale} mode, \ldots\ choose microscale burst length 
\begin{equation}
\delta\gtrsim\frac1\beta\log|\beta\Delta| .
\label{eq:deltamin}
\end{equation}
\end{itemize}

\paragraph{Use at least a second-order method}
In applications the Euler method is too inaccurate. 
We need to implement at least the Improved Euler method that is of second order accuracy in macro-step~\(\Delta\).
The Toolbox provides \verb|PIRK2()|.
Execute it without any arguments and it projectively integrates the 
multiscale Michaelis--Menten enzyme kinetics for
\(x(t)\) and~\(y(t)\),
\begin{equation*}
\frac{dx}{dt}=-x+(x+\tfrac12)y \quad\text{and}\quad
\frac{dy}{dt}=\frac1\epsilon\big[x-(x+1)y\big],
\end{equation*}
where the scale separation parameter \(\epsilon=0.05\). 
\cref{figegPIMM2} plots the results.
Here the circles, connected by dots, plot the macro-step results at time intervals \(\Delta=1\).
The solid lines are short bursts of microscale simulation used to start projecting to the next macro-time step: here the bursts are of length \(\delta=0.15\). 
Except that the first burst, by default, is twice as long to help get past the larger transients expected from far-out-of-equilibrium initial conditions.
\cref{figegPIMM2} illustrates that subsequent time-steps are in quasi-equilibrium.

Let's code a burst of length~\verb|bT| of the \ode{}s for the Michaelis--Menten enzyme kinetics at parameter~\(\epsilon\). 
First code \textsc{ode}s in a function~\verb|dMMdt| with variables \(x=\verb|x(1)|\) and \(y=\verb|x(2)|\).  
Second, starting at time~\verb|ti|, and state~\verb|xi| (row), we here simply use \Matlab's \verb|ode23| to integrate a burst \text{in time.}
\begin{matlab}
function [ts, xs] = MMburst(ti,xi,bT) 
  global MMepsilon
  dMMdt=@(t,x) [ -x(1)+(x(1)+0.5)*x(2)
    1/MMepsilon*( x(1)-(x(1)+1)*x(2))];
  [ts,xs]=ode23(dMMdt,[ti ti+bT],xi);
end
\end{matlab}
Then with initial conditions \(x(0)=1\) and \(y(0)=0\), the following script uses \verb|PIRK2()| to compute and plot a solution over time \(0\leq t\leq6\) for parameter \(\epsilon=0.05\)\,.  
Since the rate of fast decay is \(\beta\approx 1/\epsilon\) we choose a burst \text{length \(\epsilon\log(\Delta/\epsilon)\).}
\begin{matlab}
global MMepsilon
MMepsilon = 0.05
ts = 0:6
bT=MMepsilon*log((ts(2)-ts(1))/MMepsilon)
[x,tms,xms]=PIRK2(@MMburst,ts,[1;0],bT);
plot(ts,x,'o:',tms,xms)
\end{matlab}

\subsection{even integrating backward in time with forward-time simulation}

\begin{SCfigure*}
\centering
\caption{\label{figegPIMM3}projective integration \emph{backwards in time} of Michaelis--Menten enzyme kinetics from the initial condition specified at time \(t=0\).}
\includegraphics[width=0.7\textwidth, 
height=\textheight,  keepaspectratio]{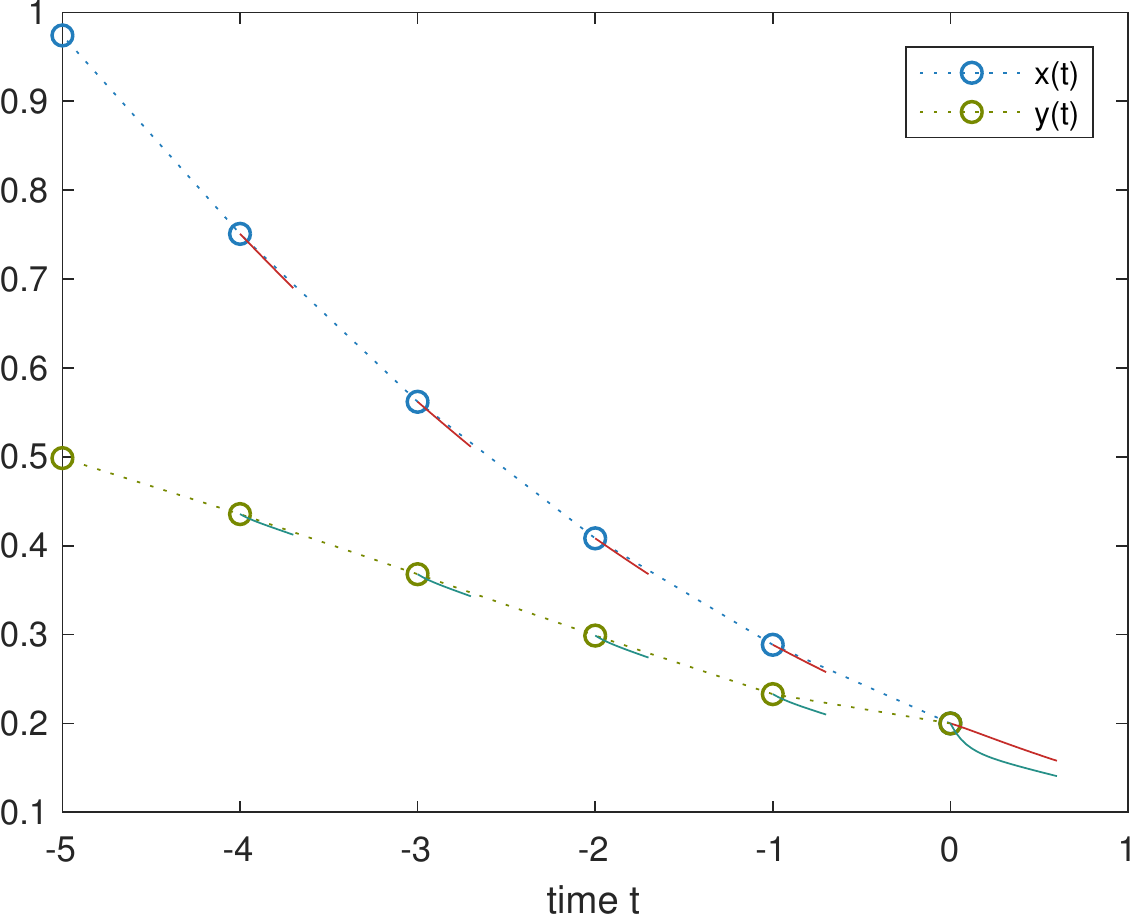}
\end{SCfigure*}%
Consider the scenario where you explore a multiscale system given by a microscopic simulator which is naturally forward in time, e.g., molecular simulator.
That is, a reverse\slash backwards simulation is not possible or not relevant (\cref{sssassb}).
Under certain conditions we may compute solutions at earlier times.
\cref{figegPIMM3} shows the example of Michaelis--Menten enzyme kinetics for which backward simulation is not feasible due to the `explosive' growth of~\(y(t)\).
Nonetheless, by simulating bursts forward in time, and then projecting backwards in time we successfully compute slow manifold solutions at earlier times.

Adapting the code of the previous subsection, the following commands draw \cref{figegPIMM3}.
The principal modification is that the macro-steps go backwards in time as in the following---also see \verb|PIRK4| with no input arguments.
\begin{matlab}
ts=0:-1:-5
bT=MMepsilon*log(abs(ts(2)-ts(1))/MMepsilon)
[xs,tms,xms]=PIRK4(@MMburst,ts,0.2*[1;1],bT);
\end{matlab}
Backward projective integration appears to need slightly longer bursts than forward, but here \(\delta= \epsilon\log(\Delta/\epsilon)\) suffices.

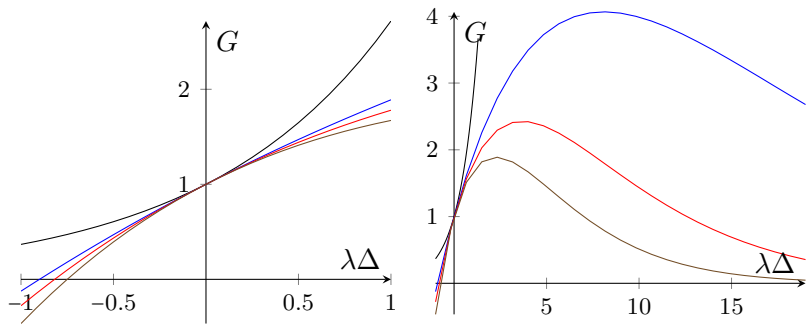
\begin{figure*}
\centering
\caption{\label{figrowthb}for scale ratio \(r=\delta/\Delta=-\tfrac19, -\tfrac29, -\tfrac39\) (blue, red, brown), plot the growth rates~\(G\) as a function of~\(\lambda\Delta\), and compare to (black)~\(e^{\lambda\Delta}\).}
\begin{tikzpicture}
\begin{axis}[small,domain=-1:1,axis lines=middle
,xlabel={$\lambda\Delta$},ylabel={$G$}]
\addplot[no marks,forget plot]{exp(x)};
\foreach \r in {-1/9,-2/9,-3/9} {
  \addplot+[no marks]{exp(\r*x)*(1+x*(1-(\r)))};
  }
\end{axis}
\end{tikzpicture}
\begin{tikzpicture}
\begin{axis}[small,domain=-1:19,axis lines=middle
,xlabel={$\lambda\Delta$},ylabel={$G$}]
\addplot[no marks,forget plot,domain=-1:1.3]{exp(x)};
\foreach \r in {-1/9,-2/9,-3/9} {
  \addplot+[no marks]{exp(\r*x)*(1+x*(1-(\r)))};
  }
\end{axis}
\end{tikzpicture}
\end{figure*}%
\cite{Gear03b} introduced the methodology.
The simple analysis of \cref{ssasss} still holds.
The difference is that here the macro-step \(\Delta<0\) and so the ratio \(r=\delta/\Delta<0\).
Consequently, in the scenario of a slow-fast multiscale system where eigenvalues~\(\lambda_i\) are either small, or large and negative, we are interested in the two cases of~\(\lambda\Delta\) small and \(\lambda\Delta\) large and positive.
\cref{figrowthb} plots the growth rate for these two cases:
\begin{itemize}
\item for macroscale accuracy choose negative time-step~\(\Delta\) such that \(\alpha|\Delta|\) is small enough;
\item to ensure the growth~\(G\) satisfies \(|G|<1\) for every microscale mode, choose microscale burst length (a large enough~\(|r|\) in \cref{figrowthb}) so that \eqref{eq:deltamin} holds (although a bit longer is better). 
\end{itemize}

\begin{activity}
Reconsider how the approximate bound~\eqref{eq:deltamin} arises from~\eqref{eq:pigrowth}.  Derive a correction to the bound for when \(|r|\)~is small but big enough to affect the right-hand side of~\eqref{eq:deltamin}.  Hence show that backward integration, negative~\(\Delta\), requires slightly longer micro-bursts than forward integration.
\end{activity}

\subsection{Projective integration via the Equation-free Toolbox}

So far we have coded three projective integration functions: \verb|PIRK2|, \verb|PIRK4|, and \verb|PIG|.

\subsubsection{Runge--Kutta-like projective integration}

The first two, as you might expect, code Runge--Kutta-like schemes of second and fourth order in the macro-step size~\(\Delta\).
They are designed to be used much like the standard \ode\ functions of \Matlab, such as \verb|ode23|.
If invoked with no arguments, then they execute the example of Michaelis--Menten enzyme kinetics: \verb|PIRK2| forwards in time; \verb|PIRK4| backwards in time.
See the example code near the start of each function.

There are some differences between \verb|PIRKn| and \Matlab\ \ode\ functions.
\begin{itemize}
\item We have not coded automatic macro-step selection, so you must specify the macro-times and steps.
\item Consequently, the times are not an output variable.
\item Instead of providing a function the computes time derivatives, you have to provide a function that computes a burst of simulation, such as \verb|MMburst| listed in \cref{ssasss}.
\item If you wish to pass the length of each burst through \verb|PIRKn| to the burst function, then supply it as the optional extra parameter. 
\item Because the burst may be of interest after the simulation, \verb|PIRKn| additionally provides optional output of the bursts.  There are three levels of burst information: none (just the macro-circles in \cref{figegPIMM2,figegPIMM3}); the physically accurate bursts (as shown in \cref{figegPIMM2,figegPIMM3}); and all computed bursts.
\end{itemize}

\begin{SCfigure*}
\centering
\caption{\label{figegPIerrs}second order errors in projective integration of Michaelis--Menten enzyme kinetics, provided the bursts are long enough.}
\includegraphics[width=0.75\textwidth, 
height=\textheight,  keepaspectratio]{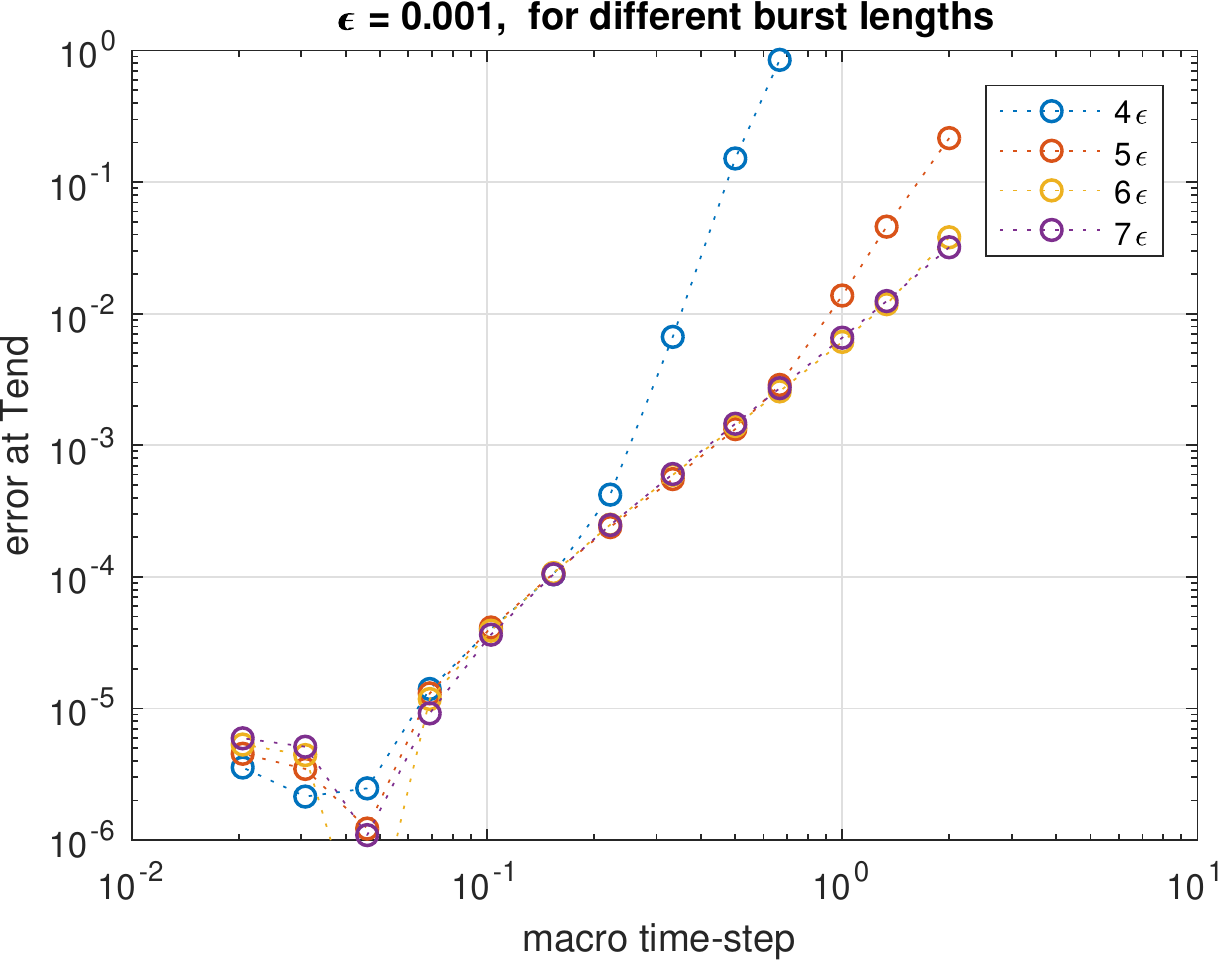}
\end{SCfigure*}

\paragraph{Errors}
Provided the microscale burst lengths are long enough, then these schemes have errors which are~\Ord{\Delta^2}, and~\Ord{\Delta^4}, correspondingly.
The script \verb|egPIerrs| illustrates this (\cref{figegPIerrs}).

For scale separation parameter \(\epsilon=0.001\), very small, \verb|egPIerrs| projectively integrates with \verb|PIRK2| the Michaelis--Menten enzyme kinetics for various macro-step lengths and burst times.
Execute \verb|egPIerrs|:
\begin{itemize}
\item for long enough bursts, the error decreases quadratically in~\(\Delta\);
\item as the burst length increases the error appears to approach the quadratic law exponentially quickly.
\end{itemize}

\subsubsection{General projective integration}

But what about adaptive codes?  Answer: we can use existing adaptive codes for macroscale integration of microscale simulations.
A user just needs a function that computes a burst of the microscale.
Then the toolbox \verb|PIG| will invoke a specified system\slash user defined function to integrate over macro-times using bursts of the microscale.

\begin{SCfigure*}
\centering
\caption{\label{figPIGsing}projective integration of a multiscale system using \texttt{ode23} on the macroscale, and \texttt{ode45} on the microscale.}
\includegraphics[width=0.7\textwidth, 
height=\textheight,  keepaspectratio]{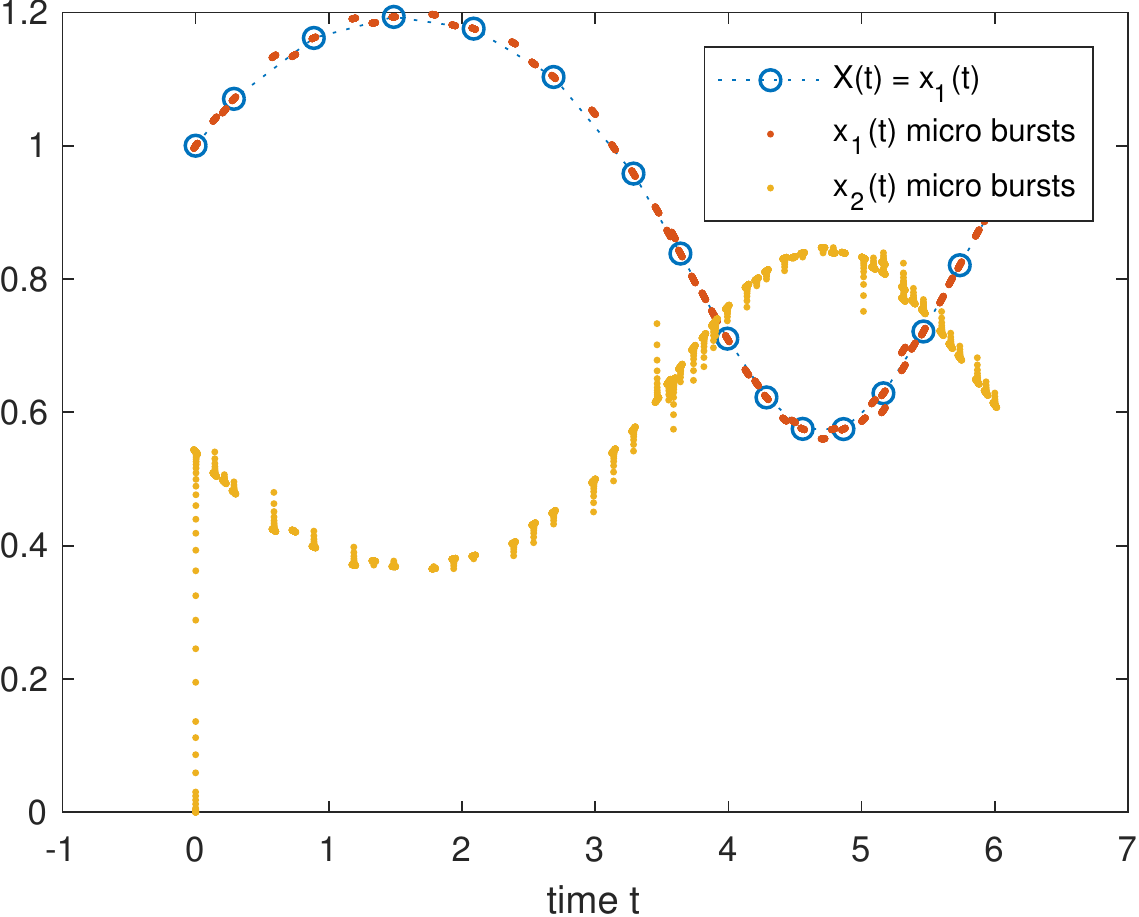}
\end{SCfigure*}

\cref{figPIGsing} shows one example (\verb|PIG| executed with no arguments):
\begin{itemize}
\item the blue circles are the macroscale computed values at macroscale times selected by the adaptive function~\verb|ode23|;
\item the red and yellow dots are the microscale bursts computed at time-steps selected by~\verb|ode45|.
\end{itemize}
The system underlying \cref{figPIGsing} is the `singular perturbation' non-autonomous \ode{}s, with parameter \(\epsilon=0.001\),
\begin{align*}&
\dot x_1=\cos x_1\sin x_2\cos t\,.
\\&
\dot x_2=\frac1\epsilon(-x_2+\cos x_1).
\end{align*}
This system is coded via (\verb|PIG|, line~214)
\begin{matlab}
epsilon = 1e-3;
dxdt=@(t,x) [ 
  cos(x(1))*sin(x(2))*cos(t)
  ( cos(x(1))-x(2) )/epsilon ];
\end{matlab}
Then here for \ode{}s the microscale burst is computed with \verb|ode45| via specifying another inline function
\begin{matlab}
bT = 2*epsilon*log(1/epsilon)
microBurst = @(tb0, xb0) feval( ...
    'ode45',dxdt,[tb0 tb0+bT],xb0);
\end{matlab}

Then invoke \verb|PIG| to execute \verb|ode23| on the coded micro-burst function over times~\([0,6]\) from an initial condition of the full microscale state.
Then plot \cref{figPIGsing}.
\begin{matlab}
x0 = [1;0];
[Ts,Xs,tms,xms] = PIG( ...
    'ode23',microBurst,[0 6] ...
    ,x0,restrict,lift);
plot(Ts,Xs,'o:',tms,xms,'.')
\end{matlab}

\paragraph{Restrict and lift between micro and macro}
The arguments \verb|restrict| and \verb|lift| invoke user specified restriction and lifting functions.  But what are they?

Notice in \cref{figPIGsing} that the macroscale (circles) is only plotted for the component~\(X(t)=x_1(t)\) and not at all for~\(x_2(t)\).
This neglect of~\(x_2\) is to show an example of restriction from the microscale to the macroscale, and a corresponding lifting from the macroscale to the microscale.
In many applications we know that macroscale quantities are relatively few in number, like temperature, pressure or patch-centre-values, whereas the microscale quantities are a morass of complexity that we have no wish to resolve over macro-times, such as molecular velocities and angular positions.
In this toy `singular perturbation' problem the natural separation is that \(x_1\)~is the slow macro-variable, and \(x_2\)~is the fast micro-variable.
Thus to \emph{restrict}\slash project microscale details into the macro-variable, that \verb|ode23| computes with, we just set \(X=\operatorname{restrict}(\xv):=x_1\).
Conversely, to \emph{lift} a macroscale state~\(X\) to a corresponding full state~\((x_1,x_2)\) we simply set \(x_1=X\) and \(x_2=x_{2,\text{approx}}\) where \(\xv_{\text{approx}}\) is some microscale state that \verb|PIG| stored from a recent micro-burst, that is, \(\xv=\operatorname{lift}(X):=(X,x_{2,\text{approx}})\): the near vertical yellow dots in \cref{figPIGsing} represent the relaxation to quasi-equilibrium from such approximations to the slow manifold.
We tell \verb|PIG| these functions via the two optional function arguments:
\begin{matlab}
restrict = @(x) x(1);
lift = @(X,xApprox) [X; xApprox(2)];
\end{matlab}

\paragraph{A methodological challenge}
Zoom in on the microscale bursts, especially the fast variable~\(x_2\): the burst looks a bit odd, it looks T-shaped.
What is going on?
Answer: the adaptive macroscale integration function expects time derivatives \emph{precisely} at the time that it specifies.  Whereas if we simulate a burst and estimate the slow derivative from the end-point of the burst then we compute a derivative at the wrong time.
So to obtain a derivative at the correct time \verb|PIG| executes two bursts: 
\footnote{Two by default, you can change.  Each of the two bursts may be shorter.}
\begin{enumerate}
\item the first burst gets to the slow manifold albeit at a wrong time; 
\item then \verb|PIG| projects backwards in time two burst-lengths; and 
\item executes a second burst which (surely) finishes at the correct time and so its estimate of the derivative is for the correct time.
\end{enumerate}

\subsection{System analysis: steady states, bifurcation, et al.}
\label{sssassb}

\paragraph{Reversing entropy}
Knock a glass of water off the table: it smashes on the floor.
Suppose we simulated with molecular dynamics.  Can we simulate backwards in time to reconstitute the glass of water?  Answer: no.
Method\slash round-off error would feed into the chaotic molecular motion so that a backwards simulation would just provide another simulation of the water spreading among the glass fragments on the floor.  
Entropy increase cannot be reversed by simulation.
\footnote{Notwithstanding \textsc{fpu} recurrence.}

Or can it?
What if we integrate a forward burst so that the macro-state variables reach the slow manifold, and then we project the macro-variables backward in time.
As in \cref{figegPIMM3} the net effect is to progress backwards in time along a slow manifold unaffected by the chaotic explosion in a direct microscale simulation. 
The nett effect is that we can integrate backwards to lower entropy states.

\begin{SCfigure*}
\centering
\caption{\label{figmacroSystem} \protect\cite[Fig.~2a]{Kevrekidis09a} projective integration of macro-state variables empowers system level analysis of chaotic microscale systems.}
\includegraphics[width=0.75\textwidth, 
height=\textheight,  keepaspectratio]{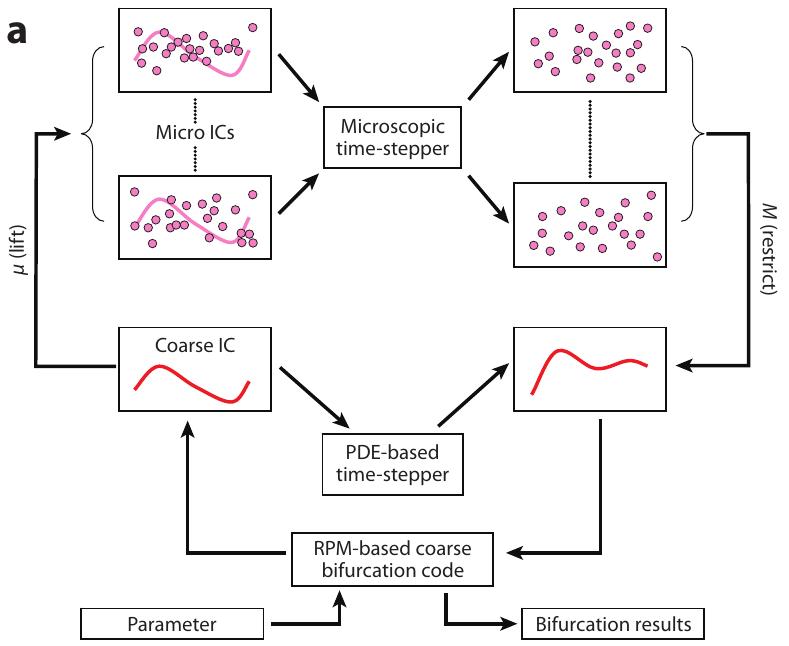}
\end{SCfigure*}%
\cite{Chiavazzo2017} extended this approach to exploring the slow manifold, rather than just one trajectory.
\cref{figmacroSystem} illustrates the idea in general: that projective integration may underpin analysis of the effective macroscale closure.
The figure also displays the possibility of parallel micro-scale simulations, from various liftings (as in ensemble simulations of weather given finite measurements), in order to better estimate subsequent macro-state variables.

\paragraph{Bifurcation analysis}
\cite{Gear02} discussed how given a macro-state, the process of lifting, a micro-burst, then restriction  (and a backward micro-time projection) results in a macro-state that we view as a map \(\Uv'=\fv(\Uv)\).
\begin{itemize}
\item Then find macro-equilibria by solving numerically \(\fv(\Uv)-\Uv=\vec 0\).
\item Determine macroscale stability, independent of microscale chaos, from the Jacobian obtained by numerically differentiating~\(\fv(\Uv)\).
\item Continuation algorithms then fill-out a bifurcation picture for the macro-variables.
\end{itemize}
One could also determine macro-state limit cycles by applying projective integration and seeking a period that repeats the macro-state.
Applications include modelling diseases \cite[]{Cisternas03}, biological dispersal \cite[]{Erban2006}, financial agents \cite[]{Siettos2012}.

Further, one could find similarity solutions by setting up equations to seek steady states of a system when space-time-state variables are scaled somehow over some small time \cite[e.g.]{Chen2004}.
For example, my script \verb|selfSim.m| uses 10,000 random walking particles to simulate diffusion on a heterogeneous period-two lattice, and then rescaled short bursts iterate to the homogenised Gaussian distribution.

\subsection{Open problems}

\begin{itemize}
\item Develop projective integration when the microscale has high frequency oscillations---homogenise over time: e.g.,  the dichotomy is \(|\lambda_i|\leq\alpha< \text{gap}<\beta\leq |\lambda_i|\).  
\item Further, develop such time-homogenisation to systems which are stochastic on micro-times---crucial.
\item Implement in the toolbox more general lifting and restricting operations.
\item But now the `baby-bathwater' question: what are appropriate macroscale variables?---after all we do not know the macroscale closure.   For example, \cite{Young01} found that Brownian Bugs should not just be modelled by densities, but also needed to model the pair correlations in order to form a qualitatively correct closure.  The challenge is to identify (dynamically?) all of the `baby' before throwing out the `bathwater'.
\item Implement `telescoping' (recursive\slash many-level) projective integration \cite[e.g.]{Gear03c}.
\item Interface the toolbox to systems analysis tools such as \textsc{auto}.
\end{itemize}

\section{Workshop: using the Equation-free Toolbox}
\localtableofcontents

This suite of \script\ functions empower users to start using the patch scheme and projective integration.
Download the current version from GitHub. 
\footnote{\protect\url{https://github.com/uoa1184615/EquationFreeGit.git}}

Many of the main functions, if invoked with \emph{no} arguments, will execute a basic example.  
For example, executing \verb|configPatches1| draws \cref{fig:burg3} arising by simulating Burgers' \pde\ within eight patches.  
Whereas executing \verb|configPatches2| computes a movie of the 2D nonlinear diffusion \(h_t=\delsq(h^3)\) on a \(9\times 7\) array of patches.

The aim of the workshop is for you to implement some example of interest to you.  Some possibilities may be inspired by examples already discussed, or the following.

\subsection{Patches in one spatial dimension}
\label{sspfeft}

The user has to drive two functions in the toolbox: \verb|configPatches1| and \verb|patchSmooth1|. 
The first helps configure the patch scheme, whereas the second provides a function to be integrated in time, or stepped in time.

We have so far designed the toolbox so that the microscale quantities are defined on a microscale lattice (although, \ldots).
In that scenario we need to create a multiscale grid in the space dimension, called~\verb|x| (a component of the struct~\verb|patches|).
Thus the \verb|configPatches1| function creates a 2D array~\verb|x| such that \(x_{ij}\)~is the \(i\)th~microscale grid point in the \(j\)th~patch.
To create this multiscale grid the user must specify: the macroscale domain,~\verb|Xlim| such as \([0,2\pi]\); the number of equi-spaced patches,~\verb|nPatch|\(=8\) say; the order of macroscale interpolation, here zero requests spectral; the (odd) number of microscale lattice points in each patch,~\verb|nSubP|\(=7\) say; and the patch micro/macro-scale \verb|ratio|, \(0.2\)~here, equal to the patch half-width divided by the \text{inter-patch spacing.}

\begin{matlab}
configPatches1(@BurgersPDE,[0 2*pi] ...
    ,nan,8,0,0.2,7);
\end{matlab}

Then a user may specify an initial condition for a simulation simply by computing an expression for all entries in \verb|patches.x|\,: for example,
\begin{matlab}
u0=0.3*(1+sin(patches.x)) ...
  +0.1*randn(size(patches.x));
\end{matlab}

During a simulation, the function \verb|patchSmooth1| computes the patch edge-values by macroscale interpolation of the patch centre-values.
The implemented order of the interpolation is specified in~\verb|ordCC|: \(2\)~is nearest neighbour quadratic; \(4\)~additionally involves the next-nearest neighbours in quartic interpolation; and so on---except that \verb|ordCC=0| is spectral interpolation.
How many edge-values are interpolated?
At least one on each edge, but if a user's microscale system is `higher-order' then it may need two or more microscale lattice edge points interpolated at each edge.
Specifying~\verb|nEdge| allows this (the default is one).

The last thing that the patch scheme needs, and the first in the parameter list for~\verb|configPatches1|, is the name of a user's function that computes microscale time derivatives\slash steps (such as the following).

\paragraph{Example of Burgers PDE inside patches}
As a microscale discretisation of Burgers' \pde\ 
\(u_t=u_{xx}-30uu_x\), here code \(\dot u_{ij} 
=\frac1{\delta x^2} (u_{i+1,j}-2u_{i,j}+u_{i-1,j}) 
-30u_{ij} \frac1{2\delta x}(u_{i+1,j}-u_{i-1,j})\).
\begin{matlab}
function ut=BurgersPDE(t,u,x)
  dx=diff(x(1:2));  % microscale spacing
  i=2:size(u,1)-1;  % inside patches
  ut=nan(size(u));  % 2D storage
  ut(i,:)=diff(u,2)/dx^2 ...
  -30*u(i,:).*(u(i+1,:)-u(i-1,:))/(2*dx);
end
\end{matlab}

\paragraph{General code overview}\ 
\begin{enumerate} 
\item invoke \verb|configPatches1()| and other initialisation
\item user time loop/integration, e.g.~\verb|ode15s|
\begin{enumerate}
\item invoke \verb|patchSmooth1()|
\begin{enumerate}
\item \verb|patchEdgeInt1()| computes edge values
\item user function for time step/derivative, e.g.~\verb|BurgersPDE()|
\end{enumerate}
\end{enumerate}
\item process results
\end{enumerate}

\begin{matlab}
[ts,ucts] = ode15s(@patchSmooth1 ...
    ,[0 0.5],u0(:));
\end{matlab}

Suppose processing the results is to draw some graphs of the simulation,
and suppose the simulation is \verb|[ts,us]=ode15s(patchSmooth1,...)|. 
The patches are most easily seen by breaking the plots between patches, as in \cref{fig:burg3}: these breaks are most easily done by assigning~\verb|nan| to the \(x\)-coordinates of the patch edges:
\verb|patches.x([1 end],:)=nan;|
\begin{itemize}
\item Then \verb|plot(patches.x(:),us(j,:)')| graphs the shape of the field at time~\(t_j\),
\item or \verb|surf(ts,patches.x(:),us')| graphs a surface over all times (as in \cref{fig:burg3}). 
\end{itemize}

\begin{activity}
Reduce the patch size by choosing a smaller ratio and see that the macroscale predictions are essentially the same.
Increase the number of points within each patch and see essentially the same.
Increase the number of patches and see the increased macroscale resolution.
\end{activity}

\begin{activity}[reaction-diffusion \pde]
Change the microscale code to solve the reaction-diffusion Ginzburg--Landau \pde\ \(u_t=u_{xx}+u-u^3\).  Use initial conditions which involve both positive and negative values of~\(u(x,0)\) and see the predicted macroscale evolves to field~\(u\) being~\(\pm1\) separated by transitions that may be relatively poorly resolved on the macroscale.
\end{activity}

\paragraph{Multiple interacting components}
Further, for systems with multiple components, such as a wave system \(du_{ij}/dt=v_{ij}\) and \(dv_{ij}/dt=\cdots\), then the microscale array~\verb|u| must be a 3D array whose third dimension has the size of the number of field variables at each microscale grid-point---generally determined from the user supplied \text{initial conditions.}

\subsection{Simulate waves on multiscale staggered grids}

\begin{SCfigure*}
\centering \caption{\label{fig:ps1WaveCtsUH}water
depth~\(h(x,t)\) (above) and velocity field~\(u(x,t)\)
(below) of the gap-tooth scheme applied to the ideal linear
wave \pde~\cref{eq:genwaveqn} with \(f_1=f_2=0\). The
microscale `random' waves persist among the propagating macroscale wave.}
\includegraphics[width=0.7\textwidth, 
height=\textheight,  keepaspectratio]{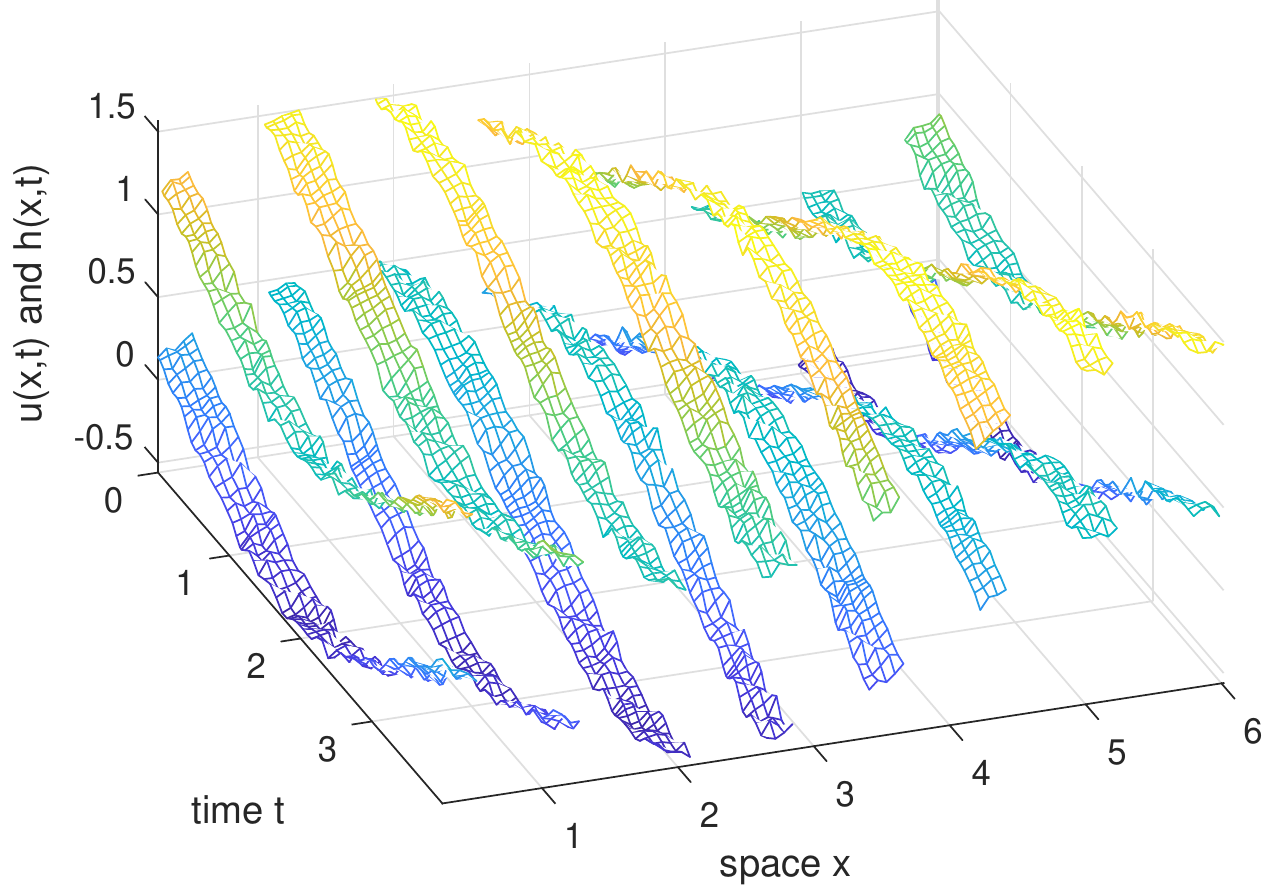}
\end{SCfigure*}

The script \verb|waterWaveExample| simulates both a linear ideal wave (\cref{fig:ps1WaveCtsUH}), and a nonlinear shallow water wave model in 1D.
To simulate the microscale detail of the waves we implement a staggered micro-grid.
Then a staggered grid of patches \cite[]{Cao2014a} empowers macroscale predictions of floods and tsunamis.

The approach developed here applies to
any wave-like system in the form
\begin{align}&
\D th=-c_1\D xu+f_1[h,u],
\nonumber\\&
\D tu=-c_2\D xh+f_2[h,u],
\label{eq:genwaveqn}
\end{align}
where the brackets indicate that the two nonlinear
functions~\(f_1\) and~\(f_2\) may involve various spatial
derivatives of the fields~\(h(x,t)\) and~\(u(x,t)\).

For wave systems, a staggered grid is best.
Let \(x_i:=i\delta x\) for microscale spacing~\(\delta x\), and grid values \(h_i(t):=h(x_i,t)\) and \(u_i(t):=u(x_i,t)\) for index~\(i\) odd/even respectively. 
Then we would code
\begin{equation*}
\begin{cases}
\dot h_i=-c_1(u_{i+1}-u_{i-1})/(2\delta x)+f_{1i} & i\text{ odd,}\\
\dot u_i=-c_2(h_{i+1}-h_{i-1})/(2\delta x)+f_{2i} & i\text{ even.}
\end{cases}
\end{equation*}
Let's implement both a staggered microscale grid and also staggered macroscale patches.
As before define \(x_{ij}:=jH+i\delta x\) for microscale spacing~\(\delta x\) and macroscale spacing~\(H\), and grid values \(h_{ij}(t):=h(x_{ij},t)\) and \(u_{ij}(t):=u(x_{ij},t)\) for index~\(i+j\) odd/even respectively. 
Then the microscale \ode{}s are
\begin{equation*}
\begin{cases}
\dot h_{ij}=-c_1(u_{i+1,j}-u_{i-1,j})/(2\delta x)+f_{1ij} & i+j\text{ odd,}\\
\dot u_{ij}=-c_2(h_{i+1,j}-h_{i-1,j})/(2\delta x)+f_{2ij} & i+j\text{ even.}
\end{cases}
\end{equation*}
With this definition, the centre-value (\(i=0\)) of the patches alternates between~\(h\) and~\(u\) values; that is, odd~\(j\) are \(h\)-patches, and even~\(j\) are \(u\)-patches.
Hence the patches are staggered.

The user's microscale code might be as in \verb|idealWavePDE.m| (a little wasteful?)
\begin{matlab}
function Ut = idealWavePDE(t,U,x)
  global patches
  dx = diff(x(2:3));
  Ut = nan(size(U));  ht = Ut;
  i = 2:size(U,1)-1;
  ht(i,:) = -(U(i+1,:)-U(i-1,:))/(2*dx);
  Ut(i,:) = -(U(i+1,:)-U(i-1,:))/(2*dx);
  Ut(patches.hPts) = ht(patches.hPts);
end
\end{matlab}

The patch\slash
gap-tooth scheme:
\begin{enumerate}
\item configPatches1, and add micro-information 
\item ode15s \(\leftrightarrow\) patchSmooth1 \(\leftrightarrow\) idealWavePDE
\item process results
\end{enumerate}
Establish the global data struct~\verb|patches| for the
\pde{}s~\cref{eq:genwaveqn} (linearised) solved on
\(2\pi\)-periodic domain, with eight patches, each patch of
half-size ratio~\(0.2\), with eleven micro-grid points
within each patch, and spectral interpolation~(\(-1\)) of
`staggered' macroscale patches to provide the edge-values of
the inter-patch coupling conditions.
\begin{matlab}
global patches
nPatch = 8
ratio = 0.2
nSubP = 11 %of the form 4*n-1
Len = 2*pi;
configPatches1(@idealWavePDE,[0 Len] ...
    ,nan,nPatch,-1,ratio,nSubP);
\end{matlab}
When the `order of interpolation' is odd, here~\(-1\), then our patch scheme interpolates the centre-values of the even patches to provide the edge-values of the odd patches, and vice-versa.

Identify and store which micro-grid points are \(h\)~or~\(u\) values
on the staggered micro-grid. 
\begin{matlab}
uPts = mod( (1:nSubP)'+(1:nPatch) ,2);
hPts = find(uPts==0);
uPts = find(uPts==1);
patches.hPts = hPts; 
patches.uPts = uPts;
\end{matlab}

Set an initial condition of some progressive wave, with noise, into~\verb|U|.
\begin{matlab}
U0 = nan(nSubP,nPatch);
U0(hPts) = 1+0.5*sin(patches.x(hPts));
U0(uPts) = 0+0.5*sin(patches.x(uPts));
U0 = U0+0.02*randn(nSubP,nPatch);
\end{matlab}

Using \verb|ode15s| we then subsample the
results because micro-grid scale waves do not dissipate and
so even \verb|ode15s| takes very small time-steps for all time---we need projective integration here.
\begin{matlab}
[ts,Ucts]=ode15s(@patchSmooth1,[0 4],U0(:));
ts = ts(1:5:end);
Ucts = Ucts(1:5:end,:);
\end{matlab}

Plot the simulation (\cref{fig:ps1WaveCtsUH}), setting \verb|nan|s to separate patches.
\begin{matlab}
xs = patches.x;  xs([1 end],:) = nan;
mesh(ts,xs(hPts),Ucts(:,hPts)'),hold on
mesh(ts,xs(uPts),Ucts(:,uPts)'),hold off
xlabel('time t'), ylabel('space x')
zlabel('u(x,t) and h(x,t)')
\end{matlab}

\begin{activity}
Reduce the patch size by choosing a smaller ratio and see that the macroscale predictions are essentially the same.
Increase the number of points within each patch and see essentially the same.
Increase the number of patches and see the increased macroscale resolution.
\end{activity}

\begin{activity}[microscale viscosity in the wave]
Change the microscale code to include viscous drag in the \(u_t=\cdots\) \pde: modify to \(u_t=-h_x+\nu u_{xx}\) for some small~\(\nu\).
Rerun and choose~\(\nu\) to see the sub-patch microscale waves effectively damped, but the macroscale waves propagating largely unaffected.
\end{activity}

\subsection{Patches of nonlinear diffusion in two space dimensions}

Similar to 1D, you have to drive two functions in the toolbox: \verb|configPatches2| and \verb|patchSmooth2|. 
The first helps configure the 2D patch scheme, whereas the second provides a function to be integrated or stepped in time.

\begin{SCfigure*}
\centering
\caption{\label{figconfigPatches2ic}a Gaussian initial condition on a \(9\times 7\) array of patches, each patch is \(5\times5\) although the edge values are not plotted.}
\includegraphics[width=0.75\textwidth]{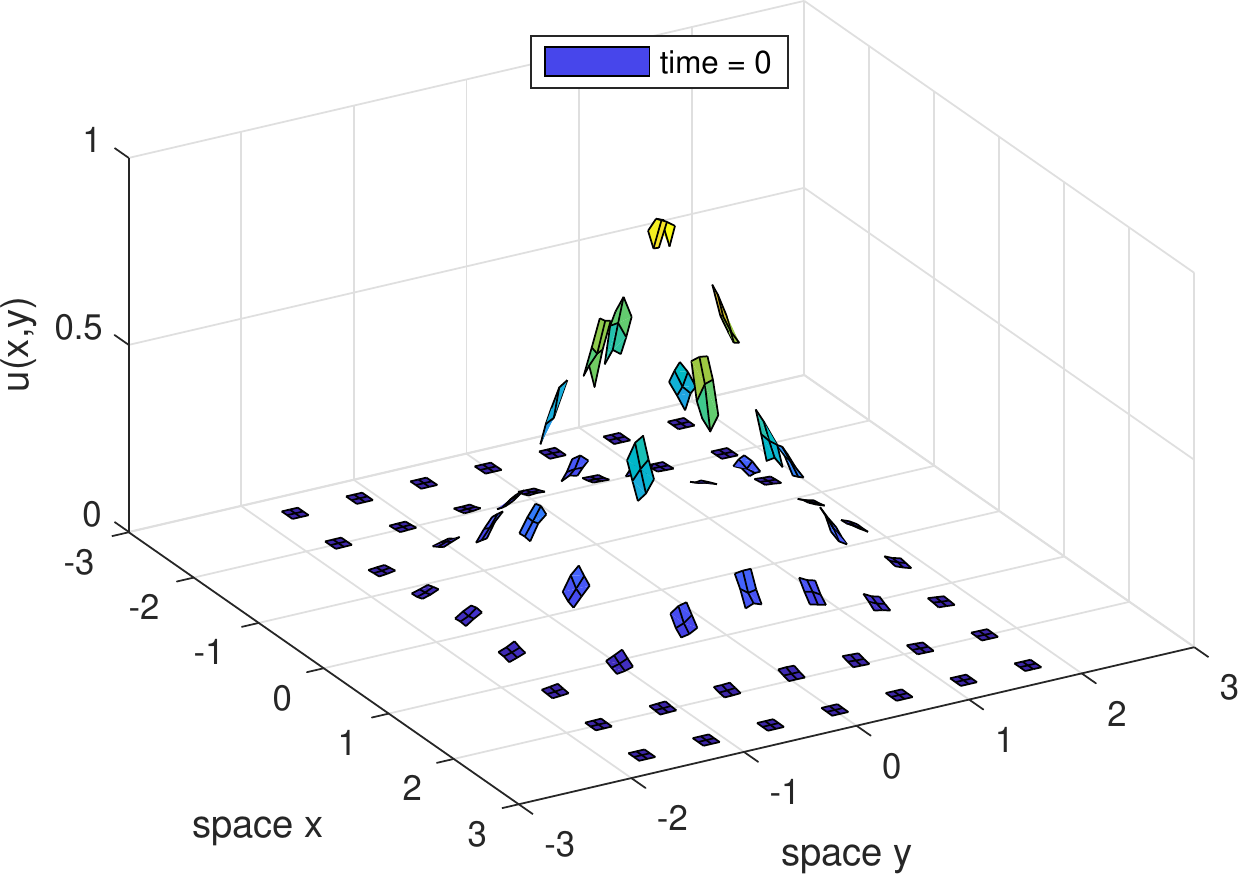}
\end{SCfigure*}%
To define 2D microscale quantities on a 2D array of patches we need to create a multiscale grid in the two space dimension (\cref{figconfigPatches2ic}), stored in~\verb|x,y| (two components of the struct~\verb|patches|).
The \verb|configPatches2| function creates two 2D arrays~\verb|x,y| such that \((x_{ik},y_{jl})\)~is the \((i,j)\)th~microscale grid point in the \((k,l)\)th~patch.
To create this multiscale grid the user must specify the 2D macroscale domain in the four elements of~\(\verb|Xlim|\), the number of equi-spaced patches in each direction in the two elements of~\verb|nPatch|, the (odd) number of microscale lattice points in each patch,~\verb|nSubP|, and the patch micro/macro-scale \verb|ratio| has two elements giving the patch half-width divided by the inter-patch spacing in each direction.
In many parameters, if the supplied parameter is a scalar, then that scalar is set for both directions. 

Then a user may specify an initial condition for a simulation by computing an expression for all entries in \verb|patches.x|\(\times\)\verb|patches.y|.
But a complication is that the microscale computation is done on a 4D array~\(u_{ijkl}\) of the field at \((x_{ik},y_{jl})\)---the \((i,j)\)th~microscale grid point in the \((k,l)\)th~patch.% 
\footnote{The reason for this subscript ordering is the thought that it should be easier to parallelise over the patches.}
So the coordinate arrays have to be rearranged into the 1st~and~3rd, and the 2nd~and~4th dimensions:
\begin{matlab}
x = reshape(patches.x,nSubP,1,[],1); 
y = reshape(patches.y,1,nSubP,1,[]);
\end{matlab}
Then auto-replication of the~\verb|x,y| arrays empowers simple assignments such as (\cref{figconfigPatches2ic})
\begin{matlab}
u0 = exp(-x.^2-y.^2);
\end{matlab}

During a simulation, the function \verb|patchSmooth2| computes the patch edge-values by macroscale interpolation of the patch centre-values.
The implemented order of the interpolation is specified in~\verb|ordCC| as in 1D, and \verb|nEdge| may specify that wider edge values are computed by the interpolation.

The last thing that the patch scheme needs, and the first in the parameter list for~\verb|configPatches2|, is the name of a user's function that computes microscale time derivatives\slash steps (such as the following).

\paragraph{Example of nonlinear diffusion PDE inside patches}
As a microscale discretisation of \(u_t=\delsq(u^3)\), code
\(\dot u_{ijkl} =\frac1{\delta x^2} (u_{i+1,j,k,l}^3
-2u_{i,j,k,l}^3 +u_{i-1,j,k,l}^3) + \frac1{\delta y^2}
(u_{i,j+1,k,l}^3 -2u_{i,j,k,l}^3 +u_{i,j-1,k,l}^3)\).
\begin{matlab}
function ut = nonDiffPDE(t,u,x,y)
dx=diff(x(1:2)); dy=diff(y(1:2)); 
i=2:size(u,1)-1; j=2:size(u,2)-1; 
ut=nan(size(u)); % preallocate
ut(i,j,:,:)=diff(u(:,j,:,:).^3,2,1)/dx^2 ...
           +diff(u(i,:,:,:).^3,2,2)/dy^2;
end
\end{matlab}

\paragraph{General code overview}
\begin{enumerate} 
\item \verb|configPatches2()| and other initialisation
\item user time loop/integration, e.g.~\verb|ode15s|
\begin{enumerate}
\item invoke \verb|patchSmooth2()|
\begin{enumerate}
\item \verb|patchEdgeInt2()| computes edge values
\item user function for time step/derivative, e.g.~\verb|nonDiffPDE|
\end{enumerate}
\end{enumerate}
\item process results
\end{enumerate}

Suppose the simulation is via \verb|[ts,us]=ode15s(patchSmooth2,...)|, and let's draw some graphs.
To graph the solution at any time~\(t_i\):
\begin{enumerate}
\item \verb|u=patchEdgeInt2(us(i,:));| converts the \(i\)th~row of~\verb|us| into a 4D array via the interpolation, 
\item then graph the macroscale, patch-centre, values with
\begin{matlab}
mesh(x((1+end)/2,1,:,1), ...
    ,y(1,(1+end)/2,1,:) ...
    ,u((1+end)/2,(1+end)/2,:,:))
\end{matlab}
\end{enumerate}
The patches are most easily seen by breaking a surface graph between patches, as in \cref{figconfigPatches2ic}:
\begin{enumerate}
\item \verb|x([1 end],:,:,:)=nan;| and \verb|y(:,[1 end],:,:)=nan;|, by assigning \verb|nan| to the \(x,y\)-coordinates of the patch edges, breaks the graphed surface;
\item \verb|u=permute(u,[1 3 2 4]);| then permutes the 4D array to separate the \(x\)-direction in the first two indices, and the \(y\)-direction in the last two; 
\item \verb|u=reshape(u,[numel(x) numel(y)]);| forms a 2D array with all the \(x\)~in the first index, and all the \(y\)~in the second;
\item \verb|surf(x(:),y(:),u')| graphs the patchy surface (as in \cref{figconfigPatches2ic}). 
\end{enumerate}

\begin{activity}[Example of 2D waves]
For \(u(x,y,t)\),  the script \verb|wave2D|, with function \verb|wavePDE|, tests and simulates the simple wave PDE in 2D space: \(u_{tt}=\delsq u\) via the two component system \(u_t=v\) and \(v_t=\delsq u\)\,.
\begin{itemize}
\item Execute and see the pure-wave nature in the pure-imaginary eigenvalues, and the wave nature in the simulation.
\item See the same macroscale simulation with smaller patches.
\item Include some dissipation in the `momentum' \pde\ (but remove the eigenvalue \verb|return| in \verb|wave2D| script: 
\begin{itemize}
\item say some simple drag \(v_t=\delsq u-0.01v\)\,;
\item then some viscous dissipation \(v_t=\delsq u +\nu\delsq v\) for small~\(\nu\approx 0.01\)\,.
\end{itemize}
\end{itemize}
\item When the viscous dissipation works for you, try an initial condition with microscale noise and see in the simulation the rapid decay of the sub-patch microscale waves.
\end{activity}

\subsection{Weave documentation in the toolbox}

To create and document the various functions, we adapt an idea due to Neil~D. Lawrence of the University of Sheffield.  
The idea is to use block comments in \Matlab\ and an environment in \LaTeX\ in order to interleave \script\ code, and its documentation in \LaTeX.  
Each function is stored in a \verb|*.m| file and has the following plan.
\begin{Verbatim}[numbers=left,xleftmargin=1.5em]
% Short explanation for "help fun"
% Author, date
%{
\section{...}
Overview LaTeX explanation.
\begin{matlab}
%}
function ...
%{
\end{matlab}
\paragraph{Input} ...
\paragraph{Output} ...
Repeated as desired:
LaTeX in end-matlab to begin-matlab
\begin{matlab}
%}
Matlab code between %} and %{
%{
\end{matlab}
Concluding LaTeX before last line.
%}
\end{Verbatim}
The function code and documentation is included by \verb|\input{*.m}| in a \LaTeX\ source file.

We need to define the environment \verb|matlab| to be some verbatim listing.  There are many available.  But \verb|fancyvrb| does a good straightforward and flexible job.
\begin{Verbatim}
\usepackage{fancyvrb}
\newenvironment{matlab}%
    {\Verbatim[numbers=left
    ,firstnumber=\the\inputlineno]}%
    {\endVerbatim}
\end{Verbatim}
Optionally, we get \verb|fancyvrb| to omit the block comment pairs \%\{~and~\%\}, 
although the following requires that the block comment pairs always be used.
\begin{Verbatim}
\makeatletter
\def\fancyvrbStartStop{%
  \edef\FancyVerbStartString
      {\@percentchar\@charrb} 
  \edef\FancyVerbStopString
      {\@percentchar\@charlb} }
\makeatother
\end{Verbatim}

\paragraph{Contributing to the toolbox}
Draft a function and example(s) as in the style of the toolbox---see the Full Developer's Manual (Appendix~B) for more detail.  Contact me.

\appendix
\section{Computer algebra codes cited in text}
\localtableofcontents
These are written in Reduce, a powerful, fast and free computer algebra package [\url{http://www.reduce-algebra.com/}].

\subsection{\texttt{homoDiff.txt}}
\label{AhomoDiff}
\begin{Verbatim}
Comment Homogenise period two diffusion. Order 40
construction finds the 1-3 case converges for
wavenumbers<2.5.  AJR, 23 May 2019;
on div; off allfac; on revpri; factor df;

maxo:=4; 
a:=am+ad; b:=am-ad;
% maxo:=20; am:=2; ad:=-1; % optional high-order case

depend uu,x,t;
let df(uu,t)=>duudt;
u1:=u2:=uu; duudt:=0;

let df(uu,x,~p)=>0 when numberp(p) and p>maxo;
for iter:=1:99 do begin
write res1:=-df(u1,t)
    +a*(-u1+for k:=0:maxo sum df(u2,x,k)/factorial(k))
    +b*(-u1+for k:=0:maxo sum df(u2,x,k)/factorial(k)*(-1)^k);
write res2:=-df(u2,t)
    +b*(-u2+for k:=0:maxo sum df(u1,x,k)/factorial(k))
    +a*(-u2+for k:=0:maxo sum df(u1,x,k)/factorial(k)*(-1)^k);

duudt:=duudt+(duudtd:=(res1+res2)/2);
u1:=u1+(u1d:=(res1-duudtd)/am/4);
u2:=u2-u1d;

if {res1,res2}={0,0} then write iter:=iter+10000;
end;

u1:=u1;
u2:=u2;
duudt:=duudt;
if maxo>10 then begin
on rounded; let df(uu,x,~p)=>z^(p/2); on list;
write coeff(duudt,z);
end;
end;
\end{Verbatim}

\subsection{\texttt{homoVibr.txt}}
\label{AhomoVibr}
\begin{Verbatim}
Comment homogenise lattice vibration problem.  For this
linear problem it appears that the resultant slow manifold
is exactly the same with just df(uu,t) changed to df(u,t,t)
in the evolution.   AJR, 25 May 2019;
on div; off allfac; on revpri; factor df,uu,vv;

maxo:=4; 
a:=am+ad; b:=am-ad;

depend vv,x,t; let df(vv,t)=>dvvdt;
depend uu,x,t; let df(uu,t)=>vv;
u1:=u2:=uu; dvvdt:=0;

let df(uu,x,~p)=>0 when numberp(p) and p>maxo;
for it:=1:9 do begin

  v1:=df(u1,t);  v2:=df(u2,t);
  write res1:=-df(v1,t)
      +a*(-u1+for k:=0:maxo sum df(u2,x,k)/factorial(k))
      +b*(-u1+for k:=0:maxo sum df(u2,x,k)/factorial(k)*(-1)^k);
  write res2:=-df(v2,t)
      +b*(-u2+for k:=0:maxo sum df(u1,x,k)/factorial(k))
      +a*(-u2+for k:=0:maxo sum df(u1,x,k)/factorial(k)*(-1)^k);
  
  dvvdt:=dvvdt+(res1+res2)/2;
  u1:=u1+(u1d:=(res1-res2)/8/am);
  u2:=u2-u1d;
  if {res1,res2}={0,0} then write it:=it+10000;
end;

u1:=u1;
u2:=u2;
dvvdt:=dvvdt;
end;
\end{Verbatim}

\subsection{\texttt{onePatchHomo.txt}}
\label{AonePatchHomo}
\begin{Verbatim}
Comment effective homogenisation closure of one patch of
two-periodic heterogeneous diffusion.  The diffusivities are
a and b, alternating.  The microscale spacing is d.  Find
that if the patch half-size n is even, then we get the
correct coefficient of the homogenisation closure.  
AJR, 27 Mar 2019 -- 27 May 2019;

procedure uint(x); 1-x^2; % for Dirchlet zero bdry at |x|=1
%procedure uint(x); 1-2/3*x-1/3*x^2; % for insulated bdry at x=-1

n:=2; % use dynamics of 2n-1 microgrid points inside the patch
nn:=2*n-1;
% Create linear operator
matrix ll(nn,nn),id(nn,nn);
for i:=1:nn do id(i,i):=1;
for i:=1:nn do ll(i,i):=-(a+b)/d^2;
for i:=1:nn-1 do ll(i,i+1):=ll(i+1,i):=
    (if evenp(i) then b else a)/d^2;
ll(1,n):=ll(1,n)+b*uint(-n*d)/d^2$
ll(nn,n):=ll(nn,n)+a*uint(n*d)/d^2$
lldsq:=ll*d^2;

% eval lambda=mu/d^2; find eqn for small one
charpoly:=factorize(det(ll*d^2-mu*id));
for j:=1:length(charpoly) do begin
    eqn:=part(charpoly,j,1);
    if sub({mu=0,d=0},eqn)=0 then write "found ",j:=j+10000;
end;
eqn:=eqn;
% two iterations to asymptotically solve and test
jac:=sub({mu=0,d=0},df(eqn,mu))$
let d^4=>0; % order of error
mu:=0$ res:=eqn$
mu:=mu-res/jac; res:=eqn;
effDiffCoeff:=mu/d^2/sub(x=0,df(uint(x),x,x)); 

end;
\end{Verbatim}

\end{document}